\documentclass[a4paper,12pt]{amsart}

\usepackage{amssymb}
\usepackage[a4paper,margin=2.5cm]{geometry}
\usepackage{mathrsfs}
\usepackage{hyperref}

\theoremstyle{plain}
\newtheorem{theorem}{Theorem}[section]
\newtheorem{corollary}[theorem]{Corollary}
\newtheorem{lemma}[theorem]{Lemma}
\newtheorem{proposition}[theorem]{Proposition}
\theoremstyle{remark}
\newtheorem{remark}[theorem]{Remark}
\numberwithin{equation}{section}

\newcommand{\R}{\mathbb R}
\newcommand{\N}{\mathbb N}
\newcommand{\Q}{\mathbb Q}
\newcommand{\C}{\mathbb C}
\newcommand{\Z}{\mathbb Z}
\newcommand{\dd}{\,\mathrm d}
\DeclareMathOperator{\Res}{Res}

\begin{document}
\title{Multivariate analogues of Ruijsenaars' hypergeometric function}
\author{Wolter Groenevelt and Rik Ledoux}
\address{Delft Institute of Applied Mathematics, Delft University of Technology, P.O. Box 5031, 2600 GA Delft, the Netherlands}
\email{w.g.m.groenevelt@tudelft.nl, rikledoux@gmail.com}

\begin{abstract}
We construct multivariate analogues of Ruijsenaars' hypergeometric
function as nested products of overlap coefficients between
eigenfunctions of twisted primitive elements in representations
of the modular double of $U_q(\mathfrak{sl}_2)$. The construction
is a hyperbolic counterpart of the nested multivariate
Askey--Wilson polynomials and functions. Using iterated
coproducts, we derive four systems of difference equations for which these multivariate functions are joint eigenfunctions. The four systems correspond to the two sets of variables and the two step sizes $iw_1$ and $iw_2$, which are interchanged by modular duality. We also prove that Ruijsenaars' one-variable transform is unitary in the parameter regime required for the construction, for positive as well as complex conjugate $w_1,w_2$.
A two-sided normalization of the kernel separates the input and output
weights, so that iterating one-variable transforms yields unitary
multivariate transforms with explicit inverse and Plancherel formulas.
\end{abstract}

\maketitle

\section{Introduction}

Ruijsenaars introduced the hyperbolic hypergeometric $R$-function in
\cite{R1} as a hyperbolic, or relativistic, analogue of the Gauss
hypergeometric function ${}_2F_1$. Its asymptotics, parameter symmetries, and
associated Hilbert space transform were subsequently developed in
\cite{R2,R3}. The $R$-function is built from the hyperbolic gamma function, which was introduced by Ruijsenaars in the context of relativistic integrable
systems \cite{RuijsenaarsADE}; up to an elementary Gaussian factor and a rescaling of the variable it coincides with Faddeev's noncompact quantum dilogarithm \cite{FKV}, see Appendix~\ref{app:quantum-dilogarithm}. The hyperbolic gamma function depends on two parameters $w_1,w_2$ with positive real parts. The $R$-function is a joint eigenfunction of four
difference operators of Askey--Wilson type and, for suitable parameters, is
the kernel of a unitary integral transform. These four operators come in two
pairs, one pair acting on each of the two variables of the $R$-function; in each pair one operator
shifts the variable by $iw_1$ and the other by $iw_2$. The $R$-function is symmetric
under interchanging $w_1$ and $w_2$, which exchanges the two shifts; this
symmetry is called modular duality. More recently, the $R$-function has been
identified with the Virasoro fusion kernel \cite{Roussillon}, and its place at
the top of a non-polynomial hyperbolic hypergeometric Askey scheme was described in \cite{LR}.

The corresponding $q$-hypergeometric theory is centered on the Askey--Wilson
polynomials and functions. Koelink and Stokman established the $L^2$-theory of
the Askey--Wilson function transform \cite{KS}, while Stokman realized
Askey--Wilson functions as formal matrix coefficients of principal-series
representations of $U_q(\mathfrak{sl}_2)$ \cite{S}. Van de Bult obtained the
hyperbolic counterpart \cite{vdB}: up to explicit normalization, Ruijsenaars'
function is an overlap coefficient between eigenfunctions of twisted primitive
elements in a representation of the modular double of
$U_q(\mathfrak{sl}_2)$ \cite{F}. The present construction builds on this overlap
interpretation and on the coalgebra structure of the modular double. A
different matrix-coefficient interpretation of a specialization of the
$R$-function was obtained in \cite{DFKSS}, where the two-particle hyperbolic
Halln\"as--Ruijsenaars wavefunction is realized using the modular-group action
in quantum Teichm\"uller theory.

There are several inequivalent ways to pass from the rank-one Askey--Wilson
and Ruijsenaars theories to several variables. The root-system approach leads
to the Weyl-group-invariant Koornwinder polynomials \cite{K} and the commuting
difference operators of van Diejen \cite{vDpoly}. On the hyperbolic side, the
Ruijsenaars--van Diejen systems provide integrable many-particle
difference operators \cite{Rcm,vDcm}; joint eigenfunctions for the
hyperbolic relativistic Calogero--Moser system were studied in
\cite{HR1,HR2,HR3}. These root-system and many-particle constructions should be
distinguished from the nested approach considered in the current paper.

In a second approach, Gasper and Rahman \cite{GR} introduced recursively
coupled products of one-variable Askey--Wilson polynomials. Iliev studied their
commuting difference operators and bispectrality \cite{I}, and Geronimo and
Iliev extended the construction to multivariable Askey--Wilson functions
\cite{GI}. Quantum-algebraic interpretations of these polynomials and
functions were given in \cite{G1,G2}; there the commuting operators arise from
the coalgebra structure of $U_q(\mathfrak{sl}_2)$. The aim of this paper is to
develop a modular-double counterpart of this nested construction, including
both its difference equations and its unitary transform
theory.

Starting from two families of twisted primitive elements, we construct their
eigenfunctions and identify the corresponding one-variable overlap coefficient $\psi_{\rho,\sigma}^{\lambda,u}(x,y)$ 
with Ruijsenaars' function. For $N\in\N$ we then define
\[
\boldsymbol\psi_{\rho,\sigma}^{\boldsymbol\lambda,\boldsymbol u}
(\boldsymbol x,\boldsymbol y)
:=
\prod_{j=1}^N
\psi_{x_{j-1},y_{j+1}}^{\lambda_j,u_j}
(x_j,y_j),
\qquad
x_0:=\rho,
\quad
y_{N+1}:=\sigma.
\]
Although this is a product of one-variable overlap coefficients, it is not a
separated product: its $j$th factor contains the neighboring variables
$x_{j-1}$ and $y_{j+1}$. This directed chain coupling is the source of
the multivariate structure. In particular, these functions are not intended as
Weyl-group-invariant eigenfunctions of a standard $N$-particle
Ruijsenaars--van Diejen system; they are modular-double counterparts of the
nested Gasper--Rahman and Geronimo--Iliev functions.

Our first main result concerns the difference equations. When
$u_1=\cdots=u_N=u$, iterated coproducts of the twisted primitive elements give
nested commuting families. Transferring their action through the overlap
pairing yields a triangular system of $N$ difference
equations, in the sense that the $j$th operator acts only in the last $j$ variables $x_{N-j+1},\ldots,x_N$; it contains up to $3^j$ shifts. Using modular
duality and reversal of the chain gives four such systems: they act in the two
sets of variables $\boldsymbol x$ and $\boldsymbol y$ and have step sizes
$iw_1$ and $iw_2$. The difference operators and their eigenvalue equations
extend by meromorphic continuation to complex $w_1,w_2$ with
$\operatorname{Re}w_1,\operatorname{Re}w_2>0$
(see Section~\ref{subsec:difference-equations-complex-periods}).

Our second main result is a unitary multivariate integral transform. We first
extend Ruijsenaars' one-variable transform \cite{R3} to the parameter domain needed for
the overlap construction, treating both positive $w_1,w_2$ and a nonreal
complex-conjugate pair $w_1,w_2$. Using a parameter symmetry of
Ruijsenaars' function, we prove that the one-variable transform
is unitary throughout the parameter domain in which we establish
its isometry (Theorem~\ref{thm:weighted-R-transform-isometry}). This holds for all positive $w_1,w_2$, and for all complex conjugate pairs. A two-sided normalization of the kernel
separates the input and output weights into factors involving only their
respective variable chains, which allows us to iterate the one-variable transform. The resulting transform is
unitary and has an explicit
inverse and Plancherel identity
(Theorem~\ref{thm:multivariate-two-sided-overlap-transform}). In contrast with
the coproduct construction of the difference equations, this transform
construction allows the parameters $u_1,\ldots,u_N$ to be independent.

The coproduct construction is also closely related to higher-rank
Askey--Wilson algebras \cite{PW,DDV}, which organize the intermediate Casimir
operators in tensor powers of $U_q(\mathfrak{sl}_2)$ and admit compatible
coproduct and braid-group actions \cite{CFPR}, and they are isomorphic to Kauffman bracket skein algebras of punctured spheres \cite{CL}.
On the integrable-model side, the multivariable Askey--Wilson polynomials of
Gasper and Rahman provide bispectral bases for modules of the $q$-Onsager
algebra, relevant for the open XXZ chain with general integrable boundary
conditions \cite{BM}. It would be natural to investigate whether the
modular-double overlap coefficients and unitary transforms constructed here
provide unitary representations of suitable higher-rank Askey--Wilson or skein
algebras, and whether the recoupling transforms diagonalize continuous-spectrum analogues of $q$-Onsager-type operators.

The structure of the paper is as follows. Section~\ref{sec:modular-double} recalls the modular double, its representations,
and the overlap pairing. Section~\ref{sec:twisted-primitive-elements}
constructs the required eigenfunctions, and
Section~\ref{sec:ruijsenaars-function} identifies their one-variable overlap
with the $R$-function and introduces the nested multivariate function.
Section~\ref{sec:askey-wilson-equations} derives its difference
equations, while Section~\ref{sec:ruijsenaars-transforms} establishes the
one-variable and multivariate transform theory. The appendix states the
properties of the hyperbolic gamma function used in the paper.

\subsection{Notation}
\label{subsec:notation}
If $\pm$ is used in the argument of a function, we mean a product, i.e.
\[
f(x\pm y)=f(x+y)f(x-y).
\]
If several $\pm$ signs occur in the same argument, they are chosen independently and all sign combinations occur, e.g.\ $G(\lambda\pm\rho\pm x)$ is the product 
\[
G(\lambda+\rho+x)G(\lambda+\rho-x)G(\lambda-\rho+x)G(\lambda-\rho-x).
\]

Throughout, $w_1,w_2$ denote the two parameters of the hyperbolic gamma function, and we set
\[
w:=\frac{w_1+w_2}{2},
\qquad
\alpha:=\frac{2\pi}{w_1w_2}.
\]
The hyperbolic gamma function is defined for $(w_1,w_2)$ in
\[
\mathcal P:=\{(w_1,w_2)\in\C^2:
\operatorname{Re}w_1>0,\ \operatorname{Re}w_2>0\}.
\]

\section{The modular double}
\label{sec:modular-double}

In this section we introduce the modular double of the quantum group $U_q(\mathfrak{sl}_2)$. The modular double was introduced by Faddeev in \cite{F} in order to combine a quantum group with its modularly dual counterpart. We use the realization and conventions from \cite{vdB}, which are adapted to the study of Ruijsenaars' hypergeometric function. 

Let $w_1,w_2\in\C$ be either a pair of positive parameters, or a pair of complex conjugate parameters with strictly positive real part and nonzero imaginary parts. These correspond to the two standard parameter regimes for representations of the modular double, $|q|=|\widetilde q|=1$ and $\bar q = \widetilde q^{-1}$, respectively; see \cite{F,F08}. Furthermore, these restrictions ensure compatibility of the representations with the contour pairing through the appropriate involution. The modular double combines two quantum groups with deformation parameters
\[
q := e^{i\pi w_1/w_2}
\qquad \text{and} \qquad 
\widetilde q := e^{i\pi w_2/w_1}
\]
into a single object. We call $\widetilde q$ the dual of the parameter $q$. In Sections~\ref{sec:modular-double}--\ref{sec:askey-wilson-equations} we assume in addition that $q^4\neq1$ and $\widetilde q^{\,4}\neq1$, i.e.\ $w_1/w_2,w_2/w_1\notin\frac12\Z$ (which is automatic for a complex conjugate pair), so that the denominators $q-q^{-1}$, $q^2-q^{-2}$ and their tilded counterparts occurring below are nonzero. This assumption is not needed in Section~\ref{sec:ruijsenaars-transforms}.

Recall that $U_q(\mathfrak{sl}_2)$ is the complex unital algebra generated by $K,\ K^{-1},\ E,\ F$,
subject to the relations
\[
KK^{-1}=K^{-1}K=1,
\]
\[
KE=q^2EK,
\qquad
KF=q^{-2}FK,\qquad 
EF-FE=\frac{K-K^{-1}}{q-q^{-1}}.
\]
The modular double of $U_q(\mathfrak{sl}_2)$ is the tensor product
	\[
	D 	:=
	U_q(\mathfrak{sl}_2)
	\otimes
	U_{\widetilde q}(\mathfrak{sl}_2)
	\]
	endowed with the standard tensor product algebra structure. Equivalently, $D$ is the algebra generated by
	\[
	K,\ K^{-1},\ E,\ F,
	\qquad
	\widetilde K,\ \widetilde K^{-1},\ \widetilde E,\ \widetilde F,
	\]
	such that the untilded generators satisfy the defining relations of
	$U_q(\mathfrak{sl}_2)$, the tilded generators satisfy the defining relations
	of $U_{\widetilde q}(\mathfrak{sl}_2)$, and the two copies commute.
We regard $D$ in this way from now on, writing $K$ for $K\otimes1$ and $\widetilde K$ for $1\otimes\widetilde K$, and similarly for the other generators. The symbol $\otimes$ will be reserved for tensor products of copies of $D$; in particular $D^{\otimes n}$ denotes the $n$-fold tensor product of the algebra $D$ with itself.

The coalgebra structure on $D$ is inherited from the coalgebra
structures on $U_q(\mathfrak{sl}_2)$ and $U_{\widetilde q}(\mathfrak{sl}_2)$.
The coproduct is the coassociative algebra homomorphism $\Delta:D\to D\otimes D$ defined on the untilded generators by
\[
\Delta(K)=K\otimes K,
\qquad
\Delta(E)=1\otimes E+E\otimes K,
\qquad
\Delta(F)=K^{-1}\otimes F+F\otimes 1,
\]
and on the tilded generators it is obtained by replacing $K,E,F$ by $\widetilde K,\widetilde E,\widetilde F$.

\subsection{Representations of the modular double}
 
We now describe the representations of the modular double on meromorphic functions used in \cite{vdB}. They are closely related to the continuous series of representations of the noncompact real form of the modular double studied in \cite{PT}.
Let $\mathcal M$ denote the space of meromorphic functions on $\C$.
For $y\in\C$, define the translation and multiplication operators
$T_y$ and $S_y$ on $\mathcal M$ by
\[
(T_y f)(x) = f(x+y),
\qquad
(S_y f)(x) = e^{2\pi ixy}f(x).
\]
For each $\lambda\in\C$, we define a representation $\pi_\lambda : D \longrightarrow \operatorname{End}(\mathcal M)$
by specifying its action on the generators of the modular double.

For the generators of $U_q(\mathfrak{sl}_2)$, we set
\begin{equation} \label{eq:pi-on-generators}
\begin{gathered}
\pi_\lambda(K)=T_{iw_1},
\qquad
\pi_\lambda(K^{-1})=T_{-iw_1}, \\
\pi_\lambda(E)
=
\frac{i q^{1/2}}{q-q^{-1}}
S_{1/(iw_2)}
\left(
q^{-1/2}e^{\pi\lambda/w_2}
+
q^{1/2}e^{-\pi\lambda/w_2}T_{iw_1}
\right),\\
\pi_\lambda(F)
=
\frac{i q^{1/2}}{q-q^{-1}}
S_{-1/(iw_2)}
\left(
q^{-1/2}e^{\pi\lambda/w_2}
+
q^{1/2}e^{-\pi\lambda/w_2}T_{-iw_1}
\right).
\end{gathered}
\end{equation}
Here $q^{1/2}:=e^{i\pi w_1/(2w_2)}$, and similarly for $\widetilde q^{1/2}$. 
The action of the tilded generators is obtained from these formulas by interchanging $w_1$ and $w_2$ and replacing $(q,K,E,F)$ by $(\widetilde q,\widetilde K,\widetilde E,\widetilde F)$.

\subsection{\texorpdfstring{$*$}{*}-structures and a sesquilinear form on \texorpdfstring{$\mathcal M$}{M}}

We now discuss the $*$-structures on the modular double $D$ and introduce
a sesquilinear form on the space $\mathcal M$ of meromorphic functions on
$\C$.

For $w_1,w_2>0$ (so that $|q|=|\widetilde q|=1$), we define a $*$-structure on  $U_q(\mathfrak{sl}_2)$ by
\[
K^* = K,
\qquad
E^* = E,
\qquad
F^* = F.
\]
This $*$-structure extends naturally to the modular double $D$ by applying it to both
the untilded and tilded copies. Here $*$ is a conjugate-linear anti-automorphism, i.e.\ $(ab)^*=b^*a^*$ and $(ca)^*=\bar c\,a^*$ for $c\in\C$; on tensor products it acts factorwise, $(a_1\otimes\cdots\otimes a_n)^*=a_1^*\otimes\cdots\otimes a_n^*$. With these conventions the coproduct is a $*$-homomorphism, $\Delta(a^*)=\Delta(a)^*$ for $a\in D$. 

For $\bar q = \widetilde q^{-1}$, the modular double admits a
second important involution $\star$, which interchanges the untilded and tilded generators. It is defined on generators by
\[
K^\star=\widetilde K,
\qquad
E^\star=\widetilde E,
\qquad
F^\star=\widetilde F,
\]
and hence $\widetilde K^\star=K, \widetilde E^\star=E, \widetilde F^\star=F$.
This involution is specific to the modular double and reflects the symmetry
between the two deformation parameters $q$ and $\widetilde q$. The involution $\star$ extends to tensor products in the same way as $*$, and the coproduct is also a $\star$-homomorphism.

Next, we define a corresponding sesquilinear form on the meromorphic functions. We start by defining some auxiliary concepts. For $f\in\mathcal M$ define its conjugate
$\overline f$ by
\[
\overline f(z)
:=
\overline{f(\overline z)}.
\]
A deformation of the real line is a contour $\mathcal C$ with an injective, piecewise smooth parametrization $ c:\R\longrightarrow\C$ such that $\operatorname{Im}c$ is bounded and $ \lim_{t\to\pm \infty}\operatorname{Re}c(t)=\pm\infty$. The contour is oriented in the direction of increasing $t$. A meromorphic function $f$ has exponential growth rate $r$ if \[ |f(x+iy)|=O(e^{r|x|}) \] as $|x|\to\infty$, locally uniformly in $y$, away from its poles. 

If $f$ and $g$ have growth rates whose sum is negative, and if $f$ and $\overline g$ have no poles on $\mathcal C$, we define 
\[
\langle f,g\rangle_{\mathcal C} := \int_{\mathcal C}f(z)\overline g(z)\,\dd z. 
\] 
The value is unchanged under admissible contour deformations that do not cross poles.
The sesquilinear form is compatible with the representation $\pi_\lambda:D\longrightarrow \operatorname{End}(\mathcal M)$ and the $*$-structures on $D$. 
\begin{lemma}
	\label{lem:sesquilinear-star-compatibility}
Let $\mathcal C$ be a deformation of the real line, and let $f,g\in\mathcal M$ be such that the sum of their exponential growth rates is smaller than $-2\pi\operatorname{Re}(1/w_2)$. Suppose that neither $f$ nor $\overline g$ has a pole in the region swept out by the contours $\mathcal C+itw_1$, $t\in[-1,1]$; in particular, there are no poles on the three contours \[ \mathcal C-iw_1, \qquad \mathcal C, \qquad \mathcal C+iw_1, \] and the shifted contours $\mathcal C\pm iw_1$ can be deformed to $\mathcal C$ without crossing a pole of $f$ or $\overline g$. Then
	\[
	\left\langle \pi_\lambda(m)f,g\right\rangle_{\mathcal C}
	=
	\begin{cases}
		\left\langle f,\pi_{\overline\lambda}(m^*)g
		\right\rangle_{\mathcal C},
		& \text{if } w_1,w_2>0, \\[0.3cm]
		\left\langle f,\pi_{\overline\lambda}(m^\star)g
		\right\rangle_{\mathcal C},
		& \text{if } \overline{w_1}=w_2,
	\end{cases}
	\]
	for all $m\in 	\operatorname{span} \{1,K,K^{-1},E,F,FK,K^{-1}E\}$.
\end{lemma}
\begin{proof} 
It suffices to consider the operators 
\[ 
T_{i\xi_1w_1}S_{\xi_2/(iw_2)}, \qquad \xi_1,\xi_2\in\{-1,0,1\}, 
\] 
since their span contains the images of the indicated elements. Changing variables and deforming the shifted contour back to $\mathcal C$ (this is where the pole condition is used; the growth condition guarantees that all integrals converge) gives 
\[ 
\left\langle T_{i\xi_1w_1}S_{\xi_2/(iw_2)}f,g \right\rangle_{\mathcal C} = \left\langle f,S_{\xi_2/(iw_2)}T_{i\xi_1w_1}g \right\rangle_{\mathcal C}. 
\] 
Using the Weyl relation between $S$ and $T$, i.e.~$S_{\xi_2/(iw_2)}T_{i\xi_1w_1}
=
e^{-2\pi i\xi_1\xi_2 w_1/w_2}
T_{i\xi_1w_1}S_{\xi_2/(iw_2)}$, and then the explicit formulas for $\pi_\lambda$, gives the $*$-adjoint relation when $w_1,w_2>0$. If $\overline{w_1}=w_2$, conjugation interchanges $w_1$ and $w_2$, which corresponds to interchanging untilded generators with the tilded ones, giving the corresponding $\star$-adjoint relation. 
\end{proof}
A similar result for $m \in \operatorname{span}\{1,\widetilde K,\widetilde K^{-1},\widetilde E,\widetilde F,\widetilde F\widetilde K,\widetilde K^{-1}\widetilde E\}$ can be obtained by interchanging $w_1$ and $w_2$; then the strip condition involves $iw_2$ instead of $iw_1$, and the growth condition becomes ``smaller than $-2\pi\operatorname{Re}(1/w_1)$''.

\begin{remark}
\label{rem:growth-condition}
The growth condition in Lemma~\ref{lem:sesquilinear-star-compatibility} is stronger than the condition needed to define $\langle f,g\rangle_{\mathcal C}$: the operators $\pi_\lambda(E)$, $\pi_\lambda(F)$ contain the multiplication operators $S_{\pm1/(iw_2)}$, which multiply by $e^{\pm2\pi x/w_2}$, a function of exponential growth rate $2\pi\operatorname{Re}(1/w_2)>0$. In all applications below, the sum of the growth rates is $-2\alpha w=-2\pi(1/w_1+1/w_2)$ (see Section~\ref{subsec:definition-psi}), which satisfies the conditions for both the untilded and the tilded elements.
\end{remark}

\begin{remark}
\label{rem:star-not-needed}
The case $\overline{w_1}=w_2$ of Lemma~\ref{lem:sesquilinear-star-compatibility}, involving the involution $\star$, will not be needed in the remainder of the paper. The difference equations are proved for positive $w_1,w_2$ using only the $*$-structure (see Lemma~\ref{lem:straight-contour}) and are then extended to complex $w_1,w_2$ by meromorphic continuation (Section~\ref{subsec:difference-equations-complex-periods}), while in Section~\ref{sec:ruijsenaars-transforms} the conjugate case is treated directly. We include the $\star$-version, together with the corresponding statements for the eigenfunctions below, for completeness, since it reflects the symmetry of the modular double under interchanging $w_1$ and $w_2$.
\end{remark}

\section{Twisted primitive elements and their eigenfunctions} \label{sec:twisted-primitive-elements}
In this section we consider eigenfunctions of twisted primitive elements. The eigenfunctions are given in terms of hyperbolic gamma functions. The hyperbolic gamma function $G=G(w_1,w_2;\cdot)$ and its zero and pole sets $Z_\pm$ are recalled in Appendix~\ref{sec:hyperbolic-gamma}. 

\subsection{Twisted primitive elements}
 Let $\rho,u,\sigma,v \in \C$. We define two elements $Y_{u,\rho}$ and $X_{v,\sigma}$ in $D$ by
\begin{align}
	Y_{u,\rho} &:= e^{2\pi  u/w_2} E+ e^{-2\pi  u/w_2} q^{-1} FK + \frac{i\nu_\rho}{q-q^{-1}}(K-1), \label{eqn:Y-def}\\
	X_{v ,\sigma} &:= e^{2\pi v /w_2} F+ e^{-2\pi v /w_2} q K^{-1} E + \frac{i\nu_\sigma}{q-q^{-1}}(K^{-1}-1),
		\label{eqn:X-def}
\end{align}
where $\nu_s = 2 \cosh \frac{2\pi s}{w_2}$. We can relate the two twisted primitive elements by introducing an involutive algebra automorphism $\theta:D \to D$, which we define on the generators by
\begin{equation} \label{eq:theta}
\theta(K)=K^{-1},
\qquad
\theta(E)=F,
\qquad
\theta(F)=E,
\end{equation}
and similarly on the tilded generators. 
Using the map $\theta$, one easily calculates from \eqref{eqn:Y-def} and \eqref{eqn:X-def} that
\[
	\theta(Y_{u,\rho}) = X_{u,\rho}.
\]
Interchanging $w_1$ and $w_2$ gives their tilded versions, i.e.
\begin{align*}
	\widetilde{Y}_{u,\rho} &:= e^{2\pi  u/w_1} \widetilde{E} + e^{-2\pi  u/w_1} \widetilde{q}^{-1} \widetilde{F}\widetilde{K} + \frac{i\widetilde{\nu}_{\rho}}{\widetilde{q}-\widetilde{q}^{-1}}(\widetilde{K}-1), \\
	\widetilde X_{v ,\sigma} &:= e^{2\pi v /w_1}  \widetilde{F}+ e^{-2\pi  v /w_1} \widetilde{q} \widetilde{K}^{-1} \widetilde{E} +\frac{i\widetilde{\nu}_{\sigma}}{\widetilde{q}-\widetilde{q}^{-1}}(\widetilde{K}^{-1}-1),
\end{align*}
where $\widetilde \nu_s := 2\cosh\frac{2\pi  s}{w_1}$. Note that $\theta(\widetilde Y_{u,\rho}) = \widetilde X_{u,\rho}$.

The elements $Y_{u,\rho}$ and $X_{v ,\sigma}$ are contained in $\operatorname{span}\{1,K, K^{-1}, E,F, FK,K^{-1} E\}$, so that we can use Lemma \ref{lem:sesquilinear-star-compatibility} to calculate the adjoints of  $\pi_\lambda(Y_{u,\rho})$ and $\pi_\lambda(X_{v,\sigma})$ with respect to the sesquilinear form $\langle \,\cdot\,,\,\cdot \rangle_{\mathcal C}$. For this we need the adjoints in $D$ with respect to the involutions $*$ and $\star$. A short calculation gives us, for positive $w_1$ and $w_2$,
\[
Y_{u,\rho}^*=Y_{\bar u ,\bar \rho}, \qquad X_{v ,\sigma}^*=X_{\bar v  ,\bar \sigma}, \qquad	\widetilde Y_{u,\rho}^*=\widetilde Y_{\bar u ,\bar \rho}, \qquad \widetilde X_{v,\sigma}^*=\widetilde X_{\bar v,\bar \sigma},
\]
and for a conjugate pair $w_1$ and $w_2$,
\[
Y_{u,\rho}^\star=\widetilde Y_{\bar u ,\bar \rho}, \qquad X_{v ,\sigma}^\star=\widetilde X_{\bar v  ,\bar \sigma}, \qquad	\widetilde Y_{u,\rho}^\star=Y_{\bar u ,\bar \rho}, \qquad \widetilde X_{v,\sigma}^\star=X_{\bar v,\bar \sigma}.
\]

\subsection{Eigenfunctions}
We determine eigenfunctions of $\pi_\lambda(Y_{u,\rho})$ and $\pi_\lambda(X_{v ,\sigma})$. We closely follow the calculation of eigenfunctions in \cite{vdB}. 
\begin{lemma}\label{lem:eigenfunctions-Y}
The function
\begin{equation*}
	H_{x,\rho}^{\lambda,u} (z) :=  \frac{G(z +u -  iw +\frac{\lambda}{2} \pm x )}{
		G(z +u -\frac{\lambda}{2} \pm \rho )
	}
\end{equation*}
satisfies
\[
\begin{split}
\pi_\lambda(Y_{u,\rho}) H_{x,\rho}^{\lambda,u} &= \mu_x^\rho H_{x,\rho}^{\lambda,u},\\ 
\pi_\lambda(\widetilde Y_{u,\rho}) H_{x,\rho}^{\lambda,u}& = \widetilde \mu_x^\rho H_{x,\rho}^{\lambda,u},
\end{split}
\]	
where
\[
\mu_x^\rho = i\frac{\nu_x-\nu_\rho}{q-q^{-1}} \qquad \text{and} \qquad \widetilde \mu_x^\rho= i\frac{\widetilde\nu_x-\widetilde\nu_\rho}{\widetilde q-\widetilde q^{-1}}.
\]
\end{lemma} 
\begin{proof}
	A direct calculation shows that the eigenvalue equation  $\pi_\lambda(Y_{u,\rho})f=\mu_{x}^\rho f$ is equivalent to 
	\[ \frac{f(z+iw_1)}{f(z)} = \frac{ \cosh\frac{\pi}{w_2} \left(z+u+\frac{\lambda}{2}-\frac{iw_2}{2}\pm x\right) }{ \cosh\frac{\pi}{w_2} \left(z+u-\frac{\lambda}{2}+\frac{iw_1}{2}\pm\rho\right) }. 
	\] 
	The defining difference equation of $G$ shows that $H_{x,\rho}^{\lambda,u}$ satisfies this relation. Interchanging $w_1$ and $w_2$ gives the tilded eigenvalue equation. 
\end{proof}

To determine eigenfunctions of $\pi_\lambda(X_{v,\sigma})$ we use the automorphism $\theta$ defined by \eqref{eq:theta}. Let $\mathcal J$ be the involutive reflection operator on $\mathcal{M}$ given by $\mathcal Jf(z) := f(-z)$.
Then
\begin{equation*}
	\mathcal JT_{iw_1}\mathcal Jf(z) = \mathcal J[z'\mapsto f(-(z'+iw_1))](z) = f(z-iw_1) = T_{-iw_1} f(z),
\end{equation*}
and
\begin{equation*}
	\mathcal JS_{1/(iw_2)}\mathcal Jf(z) = \mathcal J[z'\mapsto e^{2\pi z'/w_2} f(-z')](z) = e^{-2\pi z/w_2}f(z) = S_{-1/(iw_2)} f(z).
\end{equation*}
Using these relations, we easily check from \eqref{eq:pi-on-generators} that
\[
	\pi_\lambda(\theta(\,\cdot\,)) = \mathcal J \mathbin{\circ}  \pi_\lambda(\,\cdot\,)  \mathbin{\circ} \mathcal J
\]
holds on the generators of $U_q(\mathfrak{sl}_2)$ and, by the definition of $\theta$ on the tilded generators, also on $U_{\widetilde q}(\mathfrak{sl}_2)$. Thus it holds on $D$.
Consequently, if $f$ is an eigenfunction of $\pi_\lambda(Y)$ for some $Y\in D$, then $\mathcal Jf$ is an eigenfunction of $\pi_\lambda(\theta(Y))$ with the same eigenvalue.
Applying this observation to Lemma \ref{lem:eigenfunctions-Y} and $\theta(Y_{u,\rho})=X_{u,\rho}$ shows that $\mathcal JH^{\lambda,v}_{y,\sigma}$ is an eigenfunction of $\pi_\lambda(X_{v,\sigma})$ with eigenvalue $\mu_{y}^\sigma$, and of $\pi_\lambda(\widetilde X_{v,\sigma})$ with eigenvalue $\widetilde\mu_{y}^\sigma$.
\begin{lemma}\label{lem:eigenfunctions-F}
	The function
	\begin{equation*}
		F_{y,\sigma}^{\lambda,v}(z) := \frac{G(-z+\bar v-iw+\frac{\bar\lambda}{2}\pm \bar y)}{G(-z+\bar v-\frac{\bar\lambda}{2}\pm \bar\sigma)}
	\end{equation*}
	satisfies
	\begin{equation*}
			\begin{split}
			\pi_{\bar\lambda}(X_{v,\sigma}^*) F_{y,\sigma}^{\lambda,v} &= \mu_{\bar y}^{\bar\sigma} F_{y,\sigma}^{\lambda,v} = \overline{\mu_y^\sigma} F_{y,\sigma}^{\lambda,v},\\
			\pi_{\bar\lambda}(\widetilde X_{v,\sigma}^*) F_{y,\sigma}^{\lambda,v} &= \widetilde\mu_{\bar y}^{\bar\sigma} F_{y,\sigma}^{\lambda,v} = \overline{\widetilde\mu_y^\sigma} F_{y,\sigma}^{\lambda,v},
			\end{split}
		\end{equation*}
		in the case of positive $w_1$ and $w_2$, and 
		\begin{equation*}
			\begin{split}
			\pi_{\bar\lambda}(X_{v,\sigma}^\star) F_{y,\sigma}^{\lambda,v} &= \widetilde\mu_{\bar y}^{\bar\sigma} F_{y,\sigma}^{\lambda,v} = \overline{\mu_y^\sigma} F_{y,\sigma}^{\lambda,v},\\
			\pi_{\bar\lambda}(\widetilde X_{v,\sigma}^\star) F_{y,\sigma}^{\lambda,v} &= \mu_{\bar y}^{\bar\sigma} F_{y,\sigma}^{\lambda,v} = \overline{\widetilde\mu_y^\sigma} F_{y,\sigma}^{\lambda,v},
			\end{split}
		\end{equation*}
	
	for a conjugate pair $w_1$ and $w_2$.
\end{lemma}
\begin{proof}
	This follows from observing that $F_{y,\sigma}^{\lambda,v}=\mathcal JH_{\bar y,\bar\sigma}^{\bar\lambda,\bar v}$.
\end{proof}
For the overlap pairing introduced in the next section, we will need the conjugate meromorphic function associated with $F_{y,\sigma}^{\lambda,v}$. The symmetry \eqref{eqn:G-symmetries} gives
\begin{equation*}
	\overline{F_{y,\sigma}^{\lambda,v}}(z) = \overline{\left( \frac{G(-\bar z+\bar v-iw+\frac{\bar\lambda}{2}\pm\bar y)}{G(-\bar z+\bar v-\frac{\bar\lambda}{2}\pm\bar\sigma)}\right)}
	=\frac{G(z-v-iw-\frac{\lambda}{2}\pm y)}{G(z-v+\frac{\lambda}{2}\pm\sigma)}.
\end{equation*}

\section{Overlap coefficients and multivariate Ruijsenaars functions}
\label{sec:ruijsenaars-function}

In this section we use the eigenfunctions from the previous section to define an
overlap function. We identify this function, up to explicit normalization factors,
with Ruijsenaars' hypergeometric function. We also record the symmetries and
analytic properties that will be needed later.

\subsection{Definition and elementary symmetries}
\label{subsec:definition-psi}
We consider the sesquilinear pairing of the two eigenfunctions $H_{x,\rho}^{\lambda,u}$ and $F_{y,\sigma}^{\lambda,v}$, which may be viewed as an overlap coefficient between these two functions.
The asymptotic estimates for the hyperbolic gamma function show that
$H_{x,\rho}^{\lambda,u}$ and $F_{y,\sigma}^{\lambda,0}$ have
exponential growth rates $\alpha(\operatorname{Im}\lambda-w)$ and
$\alpha(-\operatorname{Im}\lambda-w)$, respectively.
Consequently, the product occurring in their sesquilinear pairing has
exponential growth rate $-2\alpha w$.

The poles of
$H_{x,\rho}^{\lambda,u}(z)
\overline{F_{y,\sigma}^{\lambda,0}}(z)$
that extend downwards are contained in
\begin{equation}
	Z_- -u+iw-\frac{\lambda}{2}\pm x,
	\qquad
	Z_-+iw+\frac{\lambda}{2}\pm y,
	\label{eqn:downward-poles-psi}
\end{equation}
whereas the poles that extend upwards are contained in
\begin{equation}
	Z_+-u+\frac{\lambda}{2}\pm\rho,
	\qquad
	Z_+-\frac{\lambda}{2}\pm\sigma.
	\label{eqn:upward-poles-psi}
\end{equation}
Here $Z_\pm$ are the sets defined in Appendix \ref{app:zeros-and-poles}.
Let $\mathcal C$ be a deformation of the real line separating the upward pole
sequences in \eqref{eqn:upward-poles-psi} from the downward pole sequences in
\eqref{eqn:downward-poles-psi}. We define
\begin{equation}
	\psi_{\rho,\sigma}^{\lambda,u}(x,y)
	:=
	\left\langle
	H_{x,\rho}^{\lambda,u},
	F_{y,\sigma}^{\lambda,0}
	\right\rangle_{\mathcal C},
	\label{eqn:def-psi}
\end{equation}
or
\begin{equation}
	\psi_{\rho,\sigma}^{\lambda,u}(x,y)
	=
	\int_{\mathcal C}
	\frac{G\bigl(z+u-iw+\frac{\lambda}{2}\pm x\bigr)}
	{G\bigl(z+u-\frac{\lambda}{2}\pm\rho\bigr)}
	\frac{G\bigl(z-iw-\frac{\lambda}{2}\pm y\bigr)}
	{G\bigl(z+\frac{\lambda}{2}\pm\sigma\bigr)}
	\dd z.
	\label{eqn:psi-integral}
\end{equation}
Whenever such a separating contour exists, the integral is absolutely
convergent. Its value is independent of the choice of $\mathcal C$, provided
that the contour continues to separate the same upward and downward pole
sequences. We have set $v=0$, since the corresponding overlap coefficient for general $v$ depends on $u$ and $v$ only through $u+v$.

Several symmetries follow directly from changes of variables in
\eqref{eqn:psi-integral}. The substitution $z\mapsto z-u$ gives
\begin{equation}
	\psi_{\rho,\sigma}^{\lambda,u}(x,y)
	=
	\psi_{\sigma,\rho}^{-\lambda,-u}(y,x).
	\label{eqn:psi-switch}
\end{equation}
The substitution $z\mapsto z-\frac{u}{2}+\frac{\lambda}{2}$ gives
\begin{equation}
	\psi_{\rho,\sigma}^{\lambda,u}(x,y)
	=
	\psi_{\sigma,\rho}^{u,\lambda}(x,y).
	\label{eqn:psi-lambda-u-duality}
\end{equation}
Combining \eqref{eqn:psi-switch} and \eqref{eqn:psi-lambda-u-duality} we obtain
\begin{equation} \label{eqn:psi-duality}
	\psi_{\rho,\sigma}^{\lambda,u}(x,y)
	=
	\psi_{\rho,\sigma}^{-u,-\lambda}(y,x).
\end{equation}
Finally, conjugating the integral and using the symmetries of $G$ yields
\begin{equation*}
	\overline{
		\psi_{\rho,\sigma}^{\lambda,u}(x,y)}
	=
	\psi_{\overline\sigma,\overline\rho}^{\overline\lambda,\overline u}
	(\overline y,\overline x)
	=
	\psi_{\overline\rho,\overline\sigma}^{-\overline\lambda,-\overline u}
	(\overline x,\overline y).
\end{equation*}
In particular, this implies that $\psi_{\rho,\sigma}^{\lambda,u}(x,y)$ is real-valued for $x,y,\rho,\sigma\in\R$ and $\lambda,u\in i\R$.

\subsection{Relation with Ruijsenaars' hypergeometric function}

We now relate \eqref{eqn:def-psi} to Ruijsenaars' hypergeometric
function. We use the Barnes-type integral representation introduced
in \cite[Eqns.~(1.30), (1.33)--(1.35)]{R1}, rewritten in terms of the
shifted parameters $\boldsymbol\gamma$.

Let
\[
\boldsymbol\gamma
=
(\gamma_0,\gamma_1,\gamma_2,\gamma_3)\in\C^4
\]
and define
\begin{equation}
	\widehat{\boldsymbol\gamma}
	:=
	\frac12
	\begin{pmatrix}
		1&1&1&1\\
		1&1&-1&-1\\
		1&-1&1&-1\\
		1&-1&-1&1
	\end{pmatrix}
	\boldsymbol\gamma.
	\label{eqn:gamma-duality}
\end{equation}
Throughout this paper, we use the shifted parameters $\boldsymbol\gamma$
as the parameter arguments of $R$ and $R_{\mathrm{ren}}$.
Ruijsenaars' parameters $(c_0,c_1,c_2,c_3)$ in \cite{R1} are related to
ours by $c_0=\gamma_0+w$, $c_1=\gamma_1+w_1/2$,
$c_2=\gamma_2+w_2/2$, and $c_3=\gamma_3$.
Thus $R(w_1,w_2,\boldsymbol\gamma;x,y)$ denotes Ruijsenaars'
function evaluated at these $c$-parameters; the same convention applies
to $R_{\mathrm{ren}}$.
The parameters used in the definition of $\psi$ are related to
$\boldsymbol\gamma$ by
\begin{equation}
	\boldsymbol\gamma
	=
	-i(\rho+\lambda,\rho-\lambda,\sigma-u,-\sigma-u).
	\label{eqn:gamma-parameter-correspondence}
\end{equation}
It follows that
\[
\widehat{\boldsymbol\gamma}
=
-i(\rho-u,\rho+u,\sigma+\lambda,-\sigma+\lambda).
\]
Thus $\boldsymbol\gamma\mapsto\widehat{\boldsymbol\gamma}$
corresponds to $(\lambda,u)\mapsto(-u,-\lambda)$, with
$\rho,\sigma$ unchanged. 

We write
\[
\psi_{\boldsymbol\gamma}(x,y)
:=
\psi_{\rho,\sigma}^{\lambda,u}(x,y),
\]
with the parameters related by
\eqref{eqn:gamma-parameter-correspondence}. In the present
parametrization, the integral representation of $R$ becomes
\begin{equation}
	\begin{split}
		&R\bigl(
		w_1,w_2,\boldsymbol\gamma;
		x,y
		\bigr)
		\\
		&\quad=
		\frac{1}{\sqrt{w_1w_2}}
		\frac{\prod_{j=1}^{3}
			G\bigl(i(w+\gamma_0+\gamma_j)\bigr)}
		{G(\pm x+i\gamma_0)
			G(\pm y+i\widehat\gamma_0)}
		\\
		&\qquad\quad\times
		\int_{\mathcal C}
		\frac{
			G(z+\lambda\pm x)
			G(z-u\pm y)
		}{
			G(z+iw\pm\rho)
			G(z+iw+\lambda-u\pm\sigma)
		}
		\,\dd z.
	\end{split}
	\label{eqn:R-integral-gamma-parametrization}
\end{equation}
Here $\mathcal C$ separates the upward and downward pole sequences
of the integrand.

Under the change of variables $z\mapsto
z+u-iw-\frac{\lambda}{2}$, the integral in
\eqref{eqn:R-integral-gamma-parametrization} becomes the integral
defining $\psi_{\boldsymbol\gamma}$ in
\eqref{eqn:psi-integral}. The contour is shifted by the same amount and
continues to separate the corresponding upward and downward pole sequences.
Consequently,
\begin{equation}
	\begin{split}
		R\bigl(
		w_1,w_2,\boldsymbol\gamma;
		x,y
		\bigr)
		={}&
		\frac{1}{\sqrt{w_1w_2}}
		\frac{\prod_{j=1}^{3}
			G\bigl(i(w+\gamma_0+\gamma_j)\bigr)}
		{G(\pm x+i\gamma_0)
			G(\pm y+i\widehat\gamma_0)} \,
		\psi_{\boldsymbol\gamma}(x,y).
	\end{split}
	\label{eqn:R-psi-relation-article}
\end{equation}

Ruijsenaars defines the renormalized hypergeometric function in \cite[Eq.~(2.32)]{R1}. In our parametrization, this normalization is given by \begin{equation} 
	\begin{split} R_{\mathrm{ren}}\bigl( w_1,w_2,\boldsymbol\gamma; x,y \bigr) = \frac{1}{\sqrt{w_1w_2}} \frac{ \psi_{\boldsymbol\gamma}(x,y) }{ G(\pm x+i\gamma_0) G(\pm y+i\widehat\gamma_0) }. 
	\end{split} \label{eqn:Rren-psi-relation} 
\end{equation} 
The meromorphic continuation of $R_{\mathrm{ren}}$ in its parameters and variables is established in \cite[Thm.~2.2]{R1}. Its possible poles in the variables $x$ and $y$ are contained in
 \[
\pm x=-i\gamma_j+Z_+, \qquad \pm y=-i\widehat\gamma_j+Z_+, \qquad j=0,1,2,3. 
\]
Consequently, \eqref{eqn:Rren-psi-relation} also provides a meromorphic continuation of $\psi_{\boldsymbol\gamma}$ beyond the domain on which its contour integral is directly defined. The integral representation and \eqref{eqn:Rren-psi-relation} exhibit several symmetries. In particular, $\psi_{\boldsymbol\gamma}$ and $R_{\mathrm{ren}}$ are invariant under permutations of $\gamma_1,\gamma_2,\gamma_3$ and under the interchange $w_1\leftrightarrow w_2$. Moreover, the symmetry \eqref{eqn:psi-duality} gives the duality relation
\begin{equation} 
	\psi_{\widehat{\boldsymbol\gamma}}(y,x) = \psi_{\boldsymbol\gamma}(x,y). \label{eqn:psi-gamma-duality}
\end{equation} 
The corresponding self-duality of $R$ and $R_{\mathrm{ren}}$ agrees with the self-duality obtained from the integral representation in \cite[Eq.~(2.8)]{R1}. 

A further normalization introduced by Ruijsenaars in \cite{R2} is
\begin{equation} 
	\begin{split}
		 &\mathcal E( w_1,w_2,\boldsymbol\gamma; x,y ) := \frac{\chi(\boldsymbol\gamma)}{\sqrt{w_1w_2}} \frac{ G(2x+iw)G(2y+iw) }{ G(x+i\gamma_0) G(y+i\widehat\gamma_0) }  \frac{ \psi_{\boldsymbol\gamma}(x,y) }{\prod_{j=1}^3 G(x-i\gamma_j) G(y-i\widehat\gamma_j) },
	\end{split} \label{eqn:E-psi-relation} 
\end{equation} 
where
\begin{equation} 
	\chi(\boldsymbol\gamma) := \exp\left( \frac{i\alpha}{4} \sum_{j=0}^3\gamma_j^2 - \frac{i\alpha}{8} \left(w_1^2+w_2^2+w_1w_2\right) \right). 
\end{equation}
The advantage of this normalization is its enhanced parameter symmetry. Let $W(D_4)$ act on $\boldsymbol\gamma$ by permutations and by sign changes of an even number of components. Ruijsenaars proves in \cite[Thm.~1.1]{R2} that \begin{equation} 
	\mathcal E( w_1,w_2,p\boldsymbol\gamma; x,y ) = \mathcal E( w_1,w_2,\boldsymbol\gamma; x,y ), \qquad p\in W(D_4). \label{eqn:E-D4-symmetry} 
\end{equation} 
The function $\mathcal E$ is also self-dual: 
\begin{equation*} 
	\mathcal E( w_1,w_2,\boldsymbol\gamma; x,y ) = \mathcal E( w_1,w_2,\widehat{\boldsymbol\gamma}; y,x ).
\end{equation*} 
Under the correspondence \eqref{eqn:gamma-parameter-correspondence}, the $W(D_4)$-symmetry implies, in particular, the relevant sign and parameter-pair symmetries in $\lambda,u,\rho,\sigma$. For further properties and transformations of generalized univariate hypergeometric functions related to $R_{\mathrm{ren}}$, see also \cite{BRS}.

\subsection{Multivariate overlap coefficients}

Following the quantum-algebraic construction of multivariate Askey--Wilson
polynomials and functions in \cite{G1,G2}, we construct multivariate
analogues of $\psi_{\boldsymbol\gamma}$ using iterated coproducts of the
twisted primitive elements. The construction is recursive: the spectral parameter in one
factor becomes a boundary parameter in the next factor. This coupling is the
source of the multivariate structure. 

Fix $N\in\N$. We write
\[
\boldsymbol\lambda=(\lambda_1,\ldots,\lambda_N),
\qquad
\boldsymbol u=(u_1,\ldots,u_N),
\qquad
\boldsymbol x=(x_1,\ldots,x_N),
\qquad
\boldsymbol y=(y_1,\ldots,y_N).
\]
We also fix two boundary parameters $\rho,\sigma\in\C$ and use the
conventions
\begin{equation}
	x_0:=\rho,
	\qquad
	y_{N+1}:=\sigma.
	\label{eqn:multivariate-boundary-parameters}
\end{equation}
If $u_1=\cdots=u_N=u$, we replace the bold $\boldsymbol u$ by $u$.

\subsubsection{Iterated coproducts and their eigenfunctions}

We first recall the iterated coproduct. Set $\Delta^{1}:=\Delta$, and define recursively
\begin{equation}
	\Delta^{n}
	:=
	\left(\Delta^{n-1}\otimes\operatorname{id}_D\right)
	\circ\Delta,
	\qquad n\geq2.
	\label{eqn:iterated-coproduct-definition}
\end{equation}
Thus $\Delta^n:D\longrightarrow D^{\otimes(n+1)}$. We sometimes also use $\Delta^0=\operatorname{id}_D$. 
By coassociativity, this definition is independent of the placement of
parentheses.

The elements $Y_{u,\rho}$ and $X_{0,\sigma}$ have, respectively, the right and left coideal properties
\[
\begin{split}
\Delta(Y_{u,\rho})
&=
1\otimes Y_{u,\rho}
+
Y_{u,\rho}\otimes K, \\
\Delta(X_{0,\sigma})
&=
K^{-1}\otimes X_{0,\sigma}
+
X_{0,\sigma}\otimes1.
\end{split}
\]
Then using the definition \eqref{eqn:iterated-coproduct-definition} for the first expression in each line and coassociativity, i.e.\ $\Delta^n=(\mathrm{id}\otimes\Delta^{n-1})\circ\Delta$, for the second expression, we have for $n\geq1$
\begin{equation}
	\begin{split}
	\Delta^n(Y_{u,\rho})
	&=
	1^{\otimes n}\otimes Y_{u,\rho}
	+
	\Delta^{n-1}(Y_{u,\rho})\otimes K
	=
	1\otimes\Delta^{n-1}(Y_{u,\rho})
	+
	Y_{u,\rho}\otimes K^{\otimes n}, \\
	\Delta^n(X_{0,\sigma})
	&=
	(K^{-1})^{\otimes n}\otimes X_{0,\sigma}
	+
	\Delta^{n-1}(X_{0,\sigma})\otimes1
	=
	K^{-1}\otimes\Delta^{n-1}(X_{0,\sigma})
	+
	X_{0,\sigma}\otimes1^{\otimes n}.
	\end{split}\label{eqn:iterated-coproduct}
\end{equation}
Here we used $\Delta^{n}(K)=K^{\otimes(n+1)}$.

For $j=1,\ldots,N$, we embed these iterated
coproducts into the final $j$ tensor factors of $D^{\otimes N}$ by setting
\begin{equation}
	\mathbf Y_{u,\rho}^{(j)}
	:=
	1^{\otimes(N-j)}
	\otimes
	\Delta^{j-1}(Y_{u,\rho}),
	\qquad
	\mathbf X_{\sigma}^{(j)}
	:=
	1^{\otimes(N-j)}
	\otimes
	\Delta^{j-1}(X_{0,\sigma}).
\end{equation}
Note that equation \eqref{eqn:iterated-coproduct} gives
\[
\Delta^{k-1}(X)=(\Delta^{k-j-1}\otimes\Delta^{j})(X)=(K^{-1})^{\otimes(k-j)}\otimes\Delta^{j-1}(X)+\Delta^{k-j-1}(X)\otimes1^{\otimes j},
\]
from which it readily follows that 
\begin{equation}
	\left[
	\mathbf X_{\sigma}^{(j)},
	\mathbf X_{\sigma}^{(k)}
	\right]
	=0,
	\qquad
	j,k\in\{1,\ldots,N\}.
	\label{eqn:commuting-X-family}
\end{equation}
The multivariate difference operators introduced below are obtained by
transferring the action of these represented elements through the
multivariate overlap pairing. The commutativity of the family
$\{\mathbf X_{\sigma}^{(j)}\}_{j=1}^N$ is the algebraic reason behind the
joint eigenfunction property of the multivariate overlap function; see
Remark~\ref{rem:commutativity-formal}.

Let $\mathcal M_N$ denote the space of meromorphic functions on $\C^N$, and
write $\boldsymbol z=(z_1,\ldots,z_N)$. The
representation of $D^{\otimes N}$ on $\mathcal M_N$ is
\begin{equation*}
	\pi_{\boldsymbol\lambda}
	:=
	\pi_{\lambda_1}\otimes\cdots\otimes\pi_{\lambda_N},
\end{equation*}
where the $j$th tensor factor acts on the variable $z_j$.

Define
\begin{equation}
	\mathbf H_{\boldsymbol x,\rho}^{\boldsymbol\lambda,\boldsymbol u}
	(\boldsymbol z)
	:=
	\prod_{j=1}^N
	H_{x_j,x_{j-1}}^{\lambda_j,u_j}(z_j).
	\label{eqn:multivariate-H}
\end{equation}
The nesting of its parameters is compatible with the iterated coproduct. For the following proposition, assume that $u_1=\cdots=u_N=u$.

\begin{proposition}
	\label{prop:multivariate-H-eigenfunctions}
	For $j\in\{1,\ldots,N\}$,
	\begin{equation}
		\pi_{\boldsymbol\lambda}
		\left(
		\mathbf Y_{u,x_{N-j}}^{(j)}
		\right)
		\mathbf H_{\boldsymbol x,\rho}^{\boldsymbol\lambda,u}
		=
		\mu_{x_N}^{x_{N-j}}
		\mathbf H_{\boldsymbol x,\rho}^{\boldsymbol\lambda,u}.
		\label{eqn:multivariate-H-eigenvalue}
	\end{equation}
	In particular,
	\begin{equation*}
		\pi_{\boldsymbol\lambda}
		\left(\Delta^{N-1}(Y_{u,\rho})\right)
		\mathbf H_{\boldsymbol x,\rho}^{\boldsymbol\lambda,u}
		=
		\mu_{x_N}^{\rho}
		\mathbf H_{\boldsymbol x,\rho}^{\boldsymbol\lambda,u}.
	\end{equation*}
	The same statements hold for the tilded elements, with eigenvalues
	$\widetilde\mu_{x_N}^{x_{N-j}}$ and
	$\widetilde\mu_{x_N}^{\rho}$, respectively.
\end{proposition}

\begin{proof}
	For $j=1$, equation \eqref{eqn:multivariate-H-eigenvalue} is the
	one-variable eigenvalue equation of Lemma~\ref{lem:eigenfunctions-Y} for the last factor of
	\eqref{eqn:multivariate-H}. Suppose the statement holds for some $j<N$ and put $a:=x_{N-j-1}$, $b:=x_{N-j}$. For $k\in\{1,\ldots,N\}$ let $Y_{u,a}^{[k]}$ denote $Y_{u,a}$ placed in the $k$th tensor factor, and let
	$\mathbf K^{(j)}:=1^{\otimes(N-j)}\otimes K^{\otimes j}$. By the second expression for $\Delta^j(Y_{u,a})$ in \eqref{eqn:iterated-coproduct},
	\begin{equation*}
		\mathbf Y_{u,a}^{(j+1)}
		=
		\mathbf Y_{u,a}^{(j)}
		+
		Y_{u,a}^{[N-j]}\,\mathbf K^{(j)}.
	\end{equation*}
	From the definition of $Y_{u,a}$ we obtain the identity
	\begin{equation*}
		Y_{u,a}=Y_{u,b}
		+
		\frac{i(\nu_a-\nu_b)}{q-q^{-1}}(K-1),
	\end{equation*}
	and, applying $\Delta^{j-1}$ and using $\Delta^{j-1}(K)=K^{\otimes j}$,
	$\mathbf Y_{u,a}^{(j)}=\mathbf Y_{u,b}^{(j)}+\frac{i(\nu_a-\nu_b)}{q-q^{-1}}(\mathbf K^{(j)}-1)$. Note that $\frac{i(\nu_a-\nu_b)}{q-q^{-1}}=-\mu_b^a$. Hence
	\begin{equation*}
		\mathbf Y_{u,a}^{(j+1)}
		=
		\mathbf Y_{u,b}^{(j)}
		-\mu_b^a(\mathbf K^{(j)}-1)
		+
		Y_{u,a}^{[N-j]}\,\mathbf K^{(j)}.
	\end{equation*}
	Now $\mathbf Y_{u,b}^{(j)}$ acts on the last $j$ variables with eigenvalue $\mu_{x_N}^{b}$ by the induction hypothesis, and $Y_{u,a}^{[N-j]}$ acts on the factor $H_{x_{N-j},x_{N-j-1}}^{\lambda_{N-j},u}(z_{N-j})=H_{b,a}^{\lambda_{N-j},u}(z_{N-j})$ with eigenvalue $\mu_b^a$ by Lemma~\ref{lem:eigenfunctions-Y}; since the two act on different variables and $\mathbf K^{(j)}$ acts on the last $j$ variables only, we obtain
	\begin{equation*}
		\mathbf Y_{u,a}^{(j+1)}\mathbf H_{\boldsymbol x,\rho}^{\boldsymbol\lambda,u}
		=
		\left(\mu_{x_N}^{b}+\mu_b^a\right)\mathbf H_{\boldsymbol x,\rho}^{\boldsymbol\lambda,u}
		+
		\left(-\mu_b^a+\mu_b^a\right)\mathbf K^{(j)}\mathbf H_{\boldsymbol x,\rho}^{\boldsymbol\lambda,u}
		=
		\mu_{x_N}^{a}\mathbf H_{\boldsymbol x,\rho}^{\boldsymbol\lambda,u},
	\end{equation*}
	where we used $\mu_{x_N}^{b}+\mu_b^a=\mu_{x_N}^a$. This completes the induction. The
	proof for the tilded elements is obtained by interchanging $w_1$ and $w_2$.
\end{proof}

Similarly, define
\begin{equation}
	\mathbf F_{\boldsymbol y,\sigma}^{\boldsymbol\lambda,
		\boldsymbol v}(\boldsymbol z)
	:=
	\prod_{j=1}^N F_{y_{j},y_{j+1}}^{
		\lambda_{j},v_{j}}(z_{j}).
	\label{eqn:multivariate-F}
\end{equation}
Note that the nesting is reversed. For positive
$w_1,w_2$, the functions in \eqref{eqn:multivariate-F} satisfy
\begin{equation}
	\pi_{\overline{\boldsymbol\lambda}}
	\left((\mathbf X_{\sigma}^{(j)})^*\right)
	\mathbf F_{\boldsymbol y,\sigma}^{\boldsymbol\lambda,0}
	=
	\overline{\mu_{y_{N-j+1}}^\sigma}
	\mathbf F_{\boldsymbol y,\sigma}^{\boldsymbol\lambda,0}.
	\label{eqn:multivariate-F-eigenvalue-positive}
\end{equation}
For a conjugate pair $w_1,w_2$, the corresponding identity is
\begin{equation}
	\pi_{\overline{\boldsymbol\lambda}}
	\left((\mathbf X_{\sigma}^{(j)})^\star\right)
	\mathbf F_{\boldsymbol y,\sigma}^{\boldsymbol\lambda,0}
	=
	\overline{\mu_{y_{N-j+1}}^\sigma}
	\mathbf F_{\boldsymbol y,\sigma}^{\boldsymbol\lambda,0}.
	\label{eqn:multivariate-F-eigenvalue-conjugate}
\end{equation}

\subsubsection{Definition and symmetries}

For every $j$, let $\mathcal C_j$ be a deformation of the real line that
separates the upward and downward pole sequences of the integrand defining $\psi_{x_{j-1},y_{j+1}}^{\lambda_j,u_j}
	(x_j,y_j)$. 
We define the multivariate overlap function by
\begin{equation}
	\begin{split}
		\boldsymbol\psi_{\rho,\sigma}^{\boldsymbol\lambda,
			\boldsymbol u}(\boldsymbol x,\boldsymbol y)
		:={}&
		\prod_{j=1}^N
		\left\langle
		H_{x_j,x_{j-1}}^{\lambda_j,u_j},
		F_{y_j,y_{j+1}}^{\lambda_j,0}
		\right\rangle_{\mathcal C_j}
		\\
		={}&
		\prod_{j=1}^N
		\psi_{x_{j-1},y_{j+1}}^{\lambda_j,u_j}
		(x_j,y_j).
	\end{split}
	\label{eqn:multivariate-psi}
\end{equation}
Equivalently, it is an $N$-fold integral over
$\mathcal C_1\times\cdots\times\mathcal C_N$. Although
\eqref{eqn:multivariate-psi} is a product of one-variable overlap functions,
its factors are coupled: the $j$th factor depends on the neighboring variables
$x_{j-1}$ and $y_{j+1}$. If we want to consider \eqref{eqn:multivariate-psi} as an overlap coefficient between eigenfunctions of the $\mathbf Y^{(j)}$ operators and the $(\mathbf X^{(j)})^*$ operators, we have to restrict the parameter $\boldsymbol u$ to $(u,\ldots,u)$.

The one-variable symmetry \eqref{eqn:psi-switch} gives the multivariate
symmetry
\begin{equation}
	\boldsymbol\psi_{\rho,\sigma}^{\boldsymbol\lambda,
		\boldsymbol u}(\boldsymbol x,\boldsymbol y)
	=
	\boldsymbol\psi_{\sigma,\rho}^{-\boldsymbol\lambda^{\mathrm r},
		-\boldsymbol u^{\mathrm r}}
	(\boldsymbol y^{\mathrm r},\boldsymbol x^{\mathrm r}),
	\label{eqn:multivariate-psi-reversal}
\end{equation}
where $(x_1,\ldots,x_N)^{\mathrm r}:=(x_N,\ldots,x_1)$. Indeed, applying \eqref{eqn:psi-switch} to every factor in
\eqref{eqn:multivariate-psi} and reversing the order of the product gives
\eqref{eqn:multivariate-psi-reversal}.

\section{Askey--Wilson difference operators}
\label{sec:askey-wilson-equations}

In this section we use the quantum algebra interpretation of
$\psi_{\boldsymbol\gamma}$ and an identity relating the two families $X$ and $Y$ of twisted primitive elements to derive its difference equations. This is different from the approach in \cite{vdB}, where the difference equations are obtained from the radial part of the Casimir element. After a
normalization by hyperbolic gamma functions, the difference equations become the
Askey--Wilson difference equations satisfied by Ruijsenaars' hypergeometric
function. The argument proceeds in three steps. First, we derive a three-term
relation for the eigenfunction $H_{x,\rho}^{\lambda,u}$. We then combine
this relation with the algebraic relation between the twisted primitive
elements $X_{0,\sigma}$ and $Y_{u,\rho}$, which turns the action of
$X_{0,\sigma}$ on $H_{x,\rho}^{\lambda,u}$ into a difference operator in
$x$. Finally, compatibility of the
sesquilinear form with the $*$-structure on the modular double transfers the
action of $X_{0,\sigma}$ from one entry of the pairing to the other.

Using the coproduct construction of the multivariate functions, we generalize this procedure to obtain difference equations for the multivariate Ruijsenaars functions.

\subsection{A three-term relation for the eigenfunctions}

We first express the action of $\pi_\lambda(K^{-1})$ on
$H_{x,\rho}^{\lambda,u}$ as shifts in the spectral parameter $x$.
Set
\begin{equation} \label{eqn:multivariate-A-coefficient}
	A_{x,\rho}^{\lambda}
	:=
	i\frac{
		\cosh\frac{\pi}{w_2}
		\left(x+\rho+\frac{iw_1}{2}-\lambda\right)}
	{\sinh\frac{2\pi x}{w_2}}.
\end{equation}

\begin{lemma}
	\label{lem:H-half-step}
	For $x\notin \frac{iw_2}{2}\Z$, we have
	\begin{align}
		H_{x,\rho}^{\lambda,u}\left(z-\frac{iw_1}{2}\right)
		={}&
		A_{x,\rho}^{\lambda}
		H_{x+iw_1/2,\rho+iw_1/2}^{\lambda,u}(z)
		+
		A_{-x,\rho}^{\lambda}
		H_{x-iw_1/2,\rho+iw_1/2}^{\lambda,u}(z),
		\label{eqn:H-half-step-plus}
		\\
		={}&
		A_{x,-\rho}^{\lambda}
		H_{x+iw_1/2,\rho-iw_1/2}^{\lambda,u}(z)
		+
		A_{-x,-\rho}^{\lambda}
		H_{x-iw_1/2,\rho-iw_1/2}^{\lambda,u}(z).
		\notag
	\end{align}
\end{lemma}

\begin{proof}
	Since $H_{x,\rho}^{\lambda,u}(z)
	=H_{x,\rho}^{\lambda,0}(z+u)$ and the coefficients are independent of
	$u$, it suffices to prove the first identity for $u=0$. Since
	$H_{x,\rho}^{\lambda,u}$ is even in $\rho$, the second identity then
	follows from the first. From the definition of $H$ in Lemma \ref{lem:eigenfunctions-Y} and the difference
	equation \eqref{eqn:G-difference-w1},
	\begin{equation*}
		H_{x+iw_1/2,\rho+iw_1/2}^{\lambda,0}(z)
		=
		\frac{
			\cosh\frac{\pi}{w_2}
			\left(z-iw+\frac{\lambda}{2}+x\right)}
			{\cosh\frac{\pi}{w_2}
			\left(z-\frac{\lambda}{2}+\rho\right)}
		H_{x,\rho}^{\lambda,0} \left(z-\frac{iw_1}{2}\right).
	\end{equation*}
	Using that $H_{x,\rho}^{\lambda,u}$ is even in $x$ gives the corresponding formula with the $-iw_1/2$-shift in $x$. 
	Set
	\[
	a=\frac{\pi}{w_2}
	\left(\rho+\frac{iw_1}{2}-\lambda\right),
	\qquad
	b=\frac{\pi}{w_2}
	\left(z-iw+\frac{\lambda}{2}\right),
	\qquad
	t=\frac{\pi x}{w_2}.
	\]
	Multiplying the two formulas by $A_{x,\rho}^{\lambda}$
	and $A_{-x,\rho}^{\lambda}$, respectively, and adding,
	we find that the coefficient of
	$H_{x,\rho}^{\lambda,0}(z-iw_1/2)$ is
	\[
	\frac{
		i\bigl(
		\cosh(a+t)\cosh(b+t)
		-\cosh(a-t)\cosh(b-t)
		\bigr)
	}{
		\sinh(2t)
		\cosh\frac{\pi}{w_2}
		\left(z-\frac{\lambda}{2}+\rho\right)
	}.
	\]
	The product-to-sum identity gives
	\[
	\cosh(a+t)\cosh(b+t)
	-\cosh(a-t)\cosh(b-t)
	=\sinh(a+b)\sinh(2t).
	\]
	Since $w=(w_1+w_2)/2$, we have
	\[
	a+b=
	\frac{\pi}{w_2}
	\left(z-\frac{\lambda}{2}+\rho\right)
	-\frac{i\pi}{2},
	\]
	so $i\sinh(a+b)$ equals the remaining denominator factor.
	Thus the coefficient is $1$, proving
	\eqref{eqn:H-half-step-plus}. The calculation holds away
	from poles and extends as an identity of meromorphic functions.
\end{proof}

Applying the lemma twice gives the full-step relation.

\begin{proposition} \label{prop:H-full-step}
	Away from the singularities of the coefficients,
	\begin{equation*}
		\begin{split}
			\pi_\lambda(K^{-1})H_{x,\rho}^{\lambda,u}
			={}&
			B_{x,\rho}^{\lambda}
			H_{x+iw_1,\rho}^{\lambda,u}
			+
			C_{x,\rho}^{\lambda}
			H_{x,\rho}^{\lambda,u}
			+ B_{-x,\rho}^{\lambda}
			H_{x-iw_1,\rho}^{\lambda,u},
		\end{split}
	\end{equation*}
	where
	\begin{equation}
		\begin{split}
		B_{x,\rho}^{\lambda}
		&:= A_{x,\rho}^{\lambda}
		A_{x+iw_1/2,-\rho-iw_1/2}^{\lambda} 
		=
		\frac{
			\cosh\frac{\pi}{w_2}
			\left(x\pm\rho-\lambda+\frac{iw_1}{2}\right)}
		{
			\sinh\frac{2\pi x}{w_2}
			\sinh\frac{2\pi(x+iw)}{w_2}}
		\end{split}
		\label{eqn:B-coefficient}
	\end{equation}
	and
	\begin{equation}
		\begin{split}
			C_{x,\rho}^{\lambda}
			:=&\, A_{x,\rho}^{\lambda}
			A_{-x-iw_1/2,-\rho-iw_1/2}^{\lambda}
			+
			A_{-x,\rho}^{\lambda}
			A_{x-iw_1/2,-\rho-iw_1/2}^{\lambda} \\
			=&\,
			-
			\frac{
				\cosh\frac{\pi}{w_2}
				\left(x+\rho\pm\lambda+\frac{iw_1}{2}\right)}
			{
				\sinh\frac{2\pi x}{w_2}
				\sinh\frac{2\pi(x+iw)}{w_2}}
			-\frac{
				\cosh\frac{\pi}{w_2}
				\left(x-\rho\pm\lambda-\frac{iw_1}{2}\right)}
			{
				\sinh\frac{2\pi x}{w_2}
				\sinh\frac{2\pi(x-iw)}{w_2}}.
		\end{split}
		\label{eqn:C-coefficient}
	\end{equation}
\end{proposition}

\begin{proof}
	Apply Lemma~\ref{lem:H-half-step} first to
	$H_{x,\rho}^{\lambda,u}(z-iw_1/2)$ and then to each term obtained. This
	gives the expressions of the coefficients in terms of the half-step coefficients $A_{x,\rho}^\lambda$. 
	Direct simplification gives \eqref{eqn:B-coefficient} and
	\eqref{eqn:C-coefficient}.
\end{proof}

The same statements hold with $w_1$ and $w_2$ interchanged. In addition,
the identity $F_{y,\sigma}^{\lambda,v}
=\mathcal J H_{\bar y,\bar\sigma}^{\bar\lambda,\bar v}$ and reflection give the
corresponding relation
\begin{equation*}
	\begin{split}
		\pi_{\overline\lambda}(K)F_{y,\sigma}^{\lambda,v}
		={}&
		B_{\overline y,\overline\sigma}^{\overline\lambda}
		F_{y-iw_1,\sigma}^{\lambda,v}
		+
		C_{\overline y,\overline\sigma}^{\overline\lambda}
		F_{y,\sigma}^{\lambda,v}
		+B_{-\overline y,\overline\sigma}^{\overline\lambda}
		F_{y+iw_1,\sigma}^{\lambda,v}.
	\end{split}
\end{equation*}

\subsection{A difference equation for \texorpdfstring{$\psi$}{psi} from the algebra relation}

In order to find the action of $X_{0,\sigma}$ on the eigenfunctions of $Y_{u,\rho}$, we express $X_{0,\sigma}$ in terms of $Y_{u,\rho}$ and $K^{-1}$. For $s\in\C$ set
\begin{equation}
	a_s := \frac{ q e^{2\pi s/w_2} - q^{-1}e^{-2\pi s/w_2} }{ q^2-q^{-2} }. \label{eqn:a-s}
\end{equation}
A direct computation gives
\begin{equation}
	qa_s+q^{-1}a_{-s}=e^{2\pi s/w_2},
	\qquad
	q^{-2}a_s+a_{-s}=q^{-1}e^{-2\pi s/w_2},
	\qquad
	a_s+a_{-s}=\frac{\nu_s}{q+q^{-1}}.
	\label{eqn:a-identities}
\end{equation}
For $x,\rho,\sigma\in\C$ we define
\begin{equation}
	b_s^{\rho,\sigma} := \frac{ i\nu_\rho(a_s+a_{-s})+i\nu_\sigma }{ q-q^{-1} }
	=\frac{i}{q-q^{-1}}\left(\frac{\nu_\rho\nu_s}{q+q^{-1}}+\nu_\sigma\right). \label{eqn:b-s-rho-sigma}
\end{equation}

\begin{lemma} \label{lem:relateXY}
	Let $u,\rho,\sigma\in\C$. Then
	\begin{equation}
		X_{0,\sigma} = a_{u}K^{-1}Y_{u,\rho} + a_{-u}Y_{u,\rho}K^{-1} + b_{u}^{\rho,\sigma}(K^{-1}-1). \label{eqn:relation-XY}
	\end{equation}
\end{lemma}

\begin{proof}
	From the relations $KE=q^2EK$ and $KF=q^{-2}FK$ we obtain
	\begin{align*}
		K^{-1}Y_{u,\rho}
		={}&
		q^{-2}e^{2\pi u/w_2}EK^{-1}
		+
		q e^{-2\pi u/w_2}F
		-
		\frac{i\nu_\rho}{q-q^{-1}}(K^{-1}-1),
		\\
		Y_{u,\rho}K^{-1}
		={}&
		e^{2\pi u/w_2}EK^{-1}
		+
		q^{-1}e^{-2\pi u/w_2}F
		-
		\frac{i\nu_\rho}{q-q^{-1}}(K^{-1}-1).
	\end{align*}
	By \eqref{eqn:a-identities} and \eqref{eqn:b-s-rho-sigma}, the right-hand side of \eqref{eqn:relation-XY} therefore equals
	\[
	q^{-1}EK^{-1}+F+\frac{i\nu_\sigma}{q-q^{-1}}(K^{-1}-1)
	=
	F+qK^{-1}E+\frac{i\nu_\sigma}{q-q^{-1}}(K^{-1}-1)
	=X_{0,\sigma}. \qedhere
	\]
\end{proof}
Note that $\rho$ is arbitrary in \eqref{eqn:relation-XY}: the $\rho$-dependence of $Y_{u,\rho}$ is compensated by the constant $b_u^{\rho,\sigma}$. The same computation shows that $X_{v,\sigma}$ is obtained by replacing $a_{\pm u}$ and $b_u^{\rho,\sigma}$ by $a_{\pm(u+v)}$ and $b_{u+v}^{\rho,\sigma}$; we only need the case $v=0$.

By Lemma~\ref{lem:eigenfunctions-Y} and Proposition~\ref{prop:H-full-step}, the actions of $Y_{u,\rho}$ and $K^{-1}$ on $H_{x,\rho}^{\lambda,u}$ are realized by operators acting on the spectral variable $x$. For a meromorphic function $\varphi$ of $x$ define
\begin{equation}
	\begin{split}
	(\mathsf K_\rho^\lambda\varphi)(x)
	&:=
	B_{x,\rho}^{\lambda}\varphi(x+iw_1)
	+C_{x,\rho}^{\lambda}\varphi(x)
	+B_{-x,\rho}^{\lambda}\varphi(x-iw_1),
	\\
	(\mathsf Y_\rho\varphi)(x)
	&:=
	\mu_x^\rho\varphi(x).
	\end{split}
	\label{eqn:K-Y-tau-operators}
\end{equation}
Regarding $H_{x,\rho}^{\lambda,u}(z)$ as a function of $x$ for fixed $z$, we have
$\pi_\lambda(K^{-1})H=\mathsf K_\rho^\lambda H$ and $\pi_\lambda(Y_{u,\rho})H=\mathsf Y_\rho H$. The operators $\pi_\lambda(m)$ act in $z$ with coefficients independent of $x$, so they commute with $\mathsf K_\rho^\lambda$ and $\mathsf Y_\rho$. Consequently the order of the factors is reversed:
\[
\pi_\lambda(K^{-1}Y_{u,\rho})H=\mathsf Y_\rho\mathsf K_\rho^\lambda H,
\qquad
\pi_\lambda(Y_{u,\rho}K^{-1})H=\mathsf K_\rho^\lambda\mathsf Y_\rho H.
\]
Hence Lemma~\ref{lem:relateXY} gives
\begin{equation}
	\pi_\lambda(X_{0,\sigma})H_{x,\rho}^{\lambda,u}
	=
	\mathsf D_{\rho,\sigma}^{\lambda,u}H_{x,\rho}^{\lambda,u},
	\qquad
	\mathsf D_{\rho,\sigma}^{\lambda,u}
	:=
	a_u\mathsf Y_\rho\mathsf K_\rho^\lambda
	+a_{-u}\mathsf K_\rho^\lambda\mathsf Y_\rho
	+b_u^{\rho,\sigma}\bigl(\mathsf K_\rho^\lambda-1\bigr).
	\label{eqn:D-operator}
\end{equation}

To write $\mathsf D_{\rho,\sigma}^{\lambda,u}$ explicitly, define for $\epsilon\in\{-1,0,1\}$
\begin{equation}
	\mathfrak b_{u,\epsilon}^{s,r} := \frac{ i\nu_s a_u + i\nu_{s-iw_1\epsilon}a_{-u} + i\nu_r }{ q-q^{-1} }. \label{eqn:multivariate-b-coefficient}
\end{equation}
Then $\mathfrak b_{u,0}^{s,r}=b_u^{s,r}$ and $\mathfrak b_{u,-1}^{-s,r}=\mathfrak b_{u,1}^{s,r}$, and a short computation gives
\begin{equation}
	\mathfrak b_{u,\epsilon}^{s,r}
	= \frac{2 \cosh\left( \frac{\pi}{w_2}(u+\epsilon s+r) -\frac{i\pi w_1}{2w_2} \right) \cosh\left( \frac{\pi}{w_2}(u+\epsilon s-r) -\frac{i\pi w_1}{2w_2} \right)}{\sin\left(\frac{\pi w_1}{w_2}\right)}, \qquad \epsilon=\pm1,
	\label{eqn:b-closed-form}
\end{equation}
and
\begin{equation*}
	\mathfrak b_{u,0}^{s,r} =b_u^{s,r}= \frac{1}{\sin\left(\frac{\pi w_1}{w_2}\right)} \left( \frac{ \cosh\left(\frac{2\pi s}{w_2}\right) \cosh\left(\frac{2\pi u}{w_2}\right) }{ \cos\left(\frac{\pi w_1}{w_2}\right) } + \cosh\left(\frac{2\pi r}{w_2}\right) \right).
\end{equation*}
Since $\mu_{s}^\rho=i(\nu_s-\nu_\rho)/(q-q^{-1})$, the definition \eqref{eqn:b-s-rho-sigma} gives
\begin{equation}
	a_u\mu_x^\rho+a_{-u}\mu_{x-i\epsilon w_1}^\rho+b_u^{\rho,\sigma}
	=\mathfrak b_{u,\epsilon}^{x,\sigma},
	\qquad \epsilon\in\{-1,0,1\},
	\label{eqn:mu-b-identity}
\end{equation}
which is independent of $\rho$. Set
\begin{equation*}
	\mathcal B_{\rho,\sigma}^{\lambda,u}(x)
	:=
	\mathfrak b_{u,-1}^{x,\sigma}B_{x,\rho}^{\lambda} \quad \text{and} \quad
	\mathcal C_{\rho,\sigma}^{\lambda,u}(x)
	:=
	b_u^{x,\sigma}C_{x,\rho}^{\lambda}
	-b_u^{\rho,\sigma}.
\end{equation*}
Then $\mathcal B_{\rho,\sigma}^{\lambda,u}(-x)=\mathfrak b_{u,1}^{x,\sigma}B_{-x,\rho}^{\lambda}$, and \eqref{eqn:mu-b-identity} gives
\begin{equation}
	(\mathsf D_{\rho,\sigma}^{\lambda,u}\varphi)(x)
	=
	\mathcal B_{\rho,\sigma}^{\lambda,u}(x)\varphi(x+iw_1)
	+\mathcal C_{\rho,\sigma}^{\lambda,u}(x)\varphi(x)
	+\mathcal B_{\rho,\sigma}^{\lambda,u}(-x)\varphi(x-iw_1).
	\label{eqn:D-explicit}
\end{equation}
By \eqref{eqn:B-coefficient} and \eqref{eqn:b-closed-form},
\[
\mathcal B_{\rho,\sigma}^{\lambda,u}(x)
=
\frac{2\cosh\frac{\pi}{w_2}\left(x\pm\rho-\lambda+\frac{iw_1}{2}\right)
\cosh\left(\frac{\pi}{w_2}(u-x\pm\sigma)-\frac{i\pi w_1}{2w_2}\right)}
{\sin\left(\frac{\pi w_1}{w_2}\right)\sinh\frac{2\pi x}{w_2}\sinh\frac{2\pi(x+iw)}{w_2}}.
\]

To move $\pi_\lambda(X_{0,\sigma})$ through the pairing we use Lemma~\ref{lem:sesquilinear-star-compatibility}, which requires a strip of height $2w_1$ around the contour that is free of poles. Under suitable conditions a horizontal line suffices.

\begin{lemma} \label{lem:straight-contour}
Let $w_1,w_2>0$ with $w_2>3w_1$, and put $\delta:=w-2w_1=\frac{w_2-3w_1}{2}>0$. Let $\mathcal S_\delta$ be the set of parameters $(\lambda,u,\rho,\sigma,x,y)$ with $\rho,\sigma,x,y\in\R$ and $\lambda,u\in i(-\delta,\delta)$. For parameters in $\mathcal S_\delta$ put
\[
h_-:=\max\bigl\{-\operatorname{Im}u-\tfrac12\operatorname{Im}\lambda,\ \tfrac12\operatorname{Im}\lambda\bigr\},
\qquad
h_+:=w+\min\bigl\{-\operatorname{Im}u+\tfrac12\operatorname{Im}\lambda,\ -\tfrac12\operatorname{Im}\lambda\bigr\},
\]
$c:=\frac12(h_-+h_+)$ and $\mathcal C:=\R+ic$. Then:
\begin{enumerate}
\item $f=H_{x,\rho}^{\lambda,u}$ and $g=F_{y,\sigma}^{\lambda,0}$ satisfy the hypotheses of Lemma~\ref{lem:sesquilinear-star-compatibility} with the contour $\mathcal C$;
\item for $\xi,\eta\in\{-1,0,1\}$,
\[
\langle H_{x+i\xi w_1,\rho+i\eta w_1}^{\lambda,u},F_{y,\sigma}^{\lambda,0}\rangle_{\mathcal C}
=\psi_{\rho+i\eta w_1,\sigma}^{\lambda,u}(x+i\xi w_1,y).
\]
\end{enumerate}
\end{lemma}

\begin{proof}
Since $Z_-\subset\{\operatorname{Im}z\le-w\}$ and $Z_+\subset\{\operatorname{Im}z\ge w\}$ for positive $w_1,w_2$, the downward pole sequences \eqref{eqn:downward-poles-psi} lie in $\{\operatorname{Im}z\le h_-\}$ and the upward pole sequences \eqref{eqn:upward-poles-psi} lie in $\{\operatorname{Im}z\ge h_+\}$. All poles of $f$ and of $\overline g$ belong to these sequences. The difference $h_+-h_-$ is the minimum of $w\pm\operatorname{Im}\lambda$ and $w\pm\operatorname{Im}u$, which is larger than $w-\delta=2w_1$. Hence $h_-<c-w_1$ and $h_+>c+w_1$, so the strip $\{|\operatorname{Im}z-c|\le w_1\}$, which is the region swept out by $\mathcal C+itw_1$, $t\in[-1,1]$, contains no pole of $f$ or $\overline g$. The growth condition holds by Remark~\ref{rem:growth-condition}. This proves (i).

Replacing $x$ by $x+it\xi w_1$ and $\rho$ by $\rho+it\eta w_1$, $t\in[0,1]$, moves the downward and upward pole sequences of $H$ vertically by at most $w_1$, and does not affect the poles of $\overline g$. Thus the downward poles of the integrand stay in $\{\operatorname{Im}z\le h_-+w_1\}$ and the upward poles in $\{\operatorname{Im}z\ge h_+-w_1\}$, and $\mathcal C$ separates them for all $t\in[0,1]$. Along this path the pairing is therefore the defining integral \eqref{eqn:def-psi}, which depends analytically on $t$. This proves (ii).
\end{proof}

\begin{proposition}
	\label{prop:psi-algebraic-difference}
	Let the parameters $\boldsymbol\gamma$ be given in terms of
	$(\lambda,u,\rho,\sigma)$ by
	\eqref{eqn:gamma-parameter-correspondence}. Then, as an identity of meromorphic functions,
	\[
	\mathsf D_{\rho,\sigma}^{\lambda,u}\psi_{\boldsymbol\gamma}(\,\cdot\,,y)=\mu_y^\sigma\psi_{\boldsymbol\gamma}(\,\cdot\,,y).
	\]
\end{proposition}

\begin{proof} 
	We use Lemma \ref{lem:straight-contour}. First assume $w_2>3w_1>0$ and parameters in $\mathcal S_\delta$, and let $\mathcal C$ be as in the lemma; the general case then follows by meromorphic continuation, see Proposition~\ref{prop:difference-equations-complex-periods}. By Lemma~\ref{lem:eigenfunctions-F}, Lemma~\ref{lem:sesquilinear-star-compatibility} and \eqref{eqn:D-operator},
	\[
	\mu_y^\sigma\psi_{\boldsymbol\gamma}(x,y)
	=
	\left\langle H_{x,\rho}^{\lambda,u},\pi_{\overline\lambda}(X_{0,\sigma}^*)F_{y,\sigma}^{\lambda,0}\right\rangle_{\mathcal C}
	=
	\left\langle \pi_\lambda(X_{0,\sigma})H_{x,\rho}^{\lambda,u},F_{y,\sigma}^{\lambda,0}\right\rangle_{\mathcal C}
	=
	\left\langle \mathsf D_{\rho,\sigma}^{\lambda,u}H_{x,\rho}^{\lambda,u},F_{y,\sigma}^{\lambda,0}\right\rangle_{\mathcal C}.
	\]
	By \eqref{eqn:D-explicit}, $\mathsf D_{\rho,\sigma}^{\lambda,u}H_{x,\rho}^{\lambda,u}$ is a linear combination of $H_{x+i\xi w_1,\rho}^{\lambda,u}$, $\xi\in\{-1,0,1\}$, with coefficients independent of $z$. Integrating term by term and applying Lemma~\ref{lem:straight-contour}(ii) with $\eta=0$ proves the result.
\end{proof}

\subsection{Ruijsenaars' function as an eigenfunction}

We now rewrite Proposition~\ref{prop:psi-algebraic-difference} in terms of
Ruijsenaars' hypergeometric function and obtain the fact that Ruijsenaars' function is an eigenfunction of the Askey--Wilson difference operator. 

Define
\begin{equation*}
	\mathcal A_{\boldsymbol\gamma}^{w_1,w_2}(x)
	:=
	\frac{\prod_{j=0}^3
		\cosh\frac{\pi}{w_2}
		\left(x+\tfrac{iw_1}{2}+i\gamma_j\right)}
	{\sinh\frac{2\pi x}{w_2}
		\sinh\frac{2\pi(x+iw)}{w_2}}.
\end{equation*}
The second-order Askey--Wilson difference operator acting in the $x$-variable
is
\begin{equation*}
	\begin{split}
		\mathcal L_{\boldsymbol\gamma,x}^{w_1,w_2}f(x,y)
		={}&
		\mathcal A_{\boldsymbol\gamma}^{w_1,w_2}(x)
		\bigl(f(x+iw_1,y)-f(x,y)\bigr)
		+
		\mathcal A_{\boldsymbol\gamma}^{w_1,w_2}(-x)
		\bigl(f(x-iw_1,y)-f(x,y)\bigr).
	\end{split}
\end{equation*}

Define
\begin{equation}
	\Psi_{\boldsymbol\gamma}(x,y)
	:=
	G(\pm x-i\gamma_0)
	\psi_{\boldsymbol\gamma}(x,y).
	\label{eqn:Psi-normalization}
\end{equation}
By \eqref{eqn:R-psi-relation-article}, the functions
$\Psi_{\boldsymbol\gamma}$ and
$R(w_1,w_2,\boldsymbol\gamma;x,y)$ differ by a factor
that is independent of $x$. The difference equation for $G$ gives
\begin{equation}
	\begin{split}
		\Psi_{\boldsymbol\gamma}(x+iw_1,y)
		={}&
		\frac{
			\cosh\frac{\pi}{w_2}
			\left(x+\frac{iw_1}{2}-i\gamma_0\right)}
		{
			\cosh\frac{\pi}{w_2}
			\left(x+\frac{iw_1}{2}+i\gamma_0\right)}
		G(\pm x-i\gamma_0)
		\psi_{\boldsymbol\gamma}(x+iw_1,y).
	\end{split}
	\label{eqn:Psi-shift}
\end{equation}
The following coefficient identities connect the difference equation from the
algebra relation to the Askey--Wilson operator:
\begin{equation}
	\mathcal A_{\boldsymbol\gamma}^{w_1,w_2}(x)
	\frac{
		\cosh\frac{\pi}{w_2}
		\left(x+\frac{iw_1}{2}-i\gamma_0\right)}
	{
		\cosh\frac{\pi}{w_2}
		\left(x+\frac{iw_1}{2}+i\gamma_0\right)}
	=
	-\frac{i(q-q^{-1})}{4}
	\mathcal B_{\rho,\sigma}^{\lambda,u}(x)
	\label{eqn:AW-B-identity}
\end{equation}
and
\begin{equation}
	\begin{split}
		-\mathcal A_{\boldsymbol\gamma}^{w_1,w_2}(x)
		-\mathcal A_{\boldsymbol\gamma}^{w_1,w_2}(-x)
		=
		-\frac{i(q-q^{-1})}{4}
		\mathcal C_{\rho,\sigma}^{\lambda,u}(x)
		+
		\cosh\frac{\pi}{w_2}
		\left(\rho+\frac{iw_1}{2}-u\pm\sigma\right).
	\end{split}
	\label{eqn:AW-C-identity}
\end{equation}
Both identities follow by inserting the parameter correspondence
\eqref{eqn:gamma-parameter-correspondence} and applying elementary product
formulas for hyperbolic functions.

We are now ready to give a proof of the eigenvalue equations for the $R$-function based on Proposition \ref{prop:psi-algebraic-difference}. Define
\begin{equation*}
	v(w_1,w_2,\boldsymbol\gamma;y)
	:=
	\frac12\left(
	\cosh\frac{2\pi y}{w_2}
	+
	\cosh\frac{\pi i}{w_2}
	\left(w_1+2\widehat\gamma_0\right)
	\right),
\end{equation*}
where $\widehat{\boldsymbol\gamma}$ is defined by \eqref{eqn:gamma-duality}.
\begin{theorem}[Ruijsenaars \cite{R1}]
\label{thm:R-askey-wilson}
Ruijsenaars' hypergeometric function satisfies
\begin{equation}
\begin{split}
\mathcal L_{\boldsymbol\gamma,x}^{w_1,w_2}
R\bigl(w_1,w_2,\boldsymbol\gamma;x,y\bigr)
&=
v(w_1,w_2,\boldsymbol\gamma;y)
R\bigl(w_1,w_2,\boldsymbol\gamma;x,y\bigr),
 \\
\mathcal L_{\widehat{\boldsymbol\gamma},y}^{w_1,w_2}
R\bigl(w_1,w_2,\boldsymbol\gamma;x,y\bigr)
&=
v(w_1,w_2,\widehat{\boldsymbol\gamma};x)
R\bigl(w_1,w_2,\boldsymbol\gamma;x,y\bigr),
\end{split}
\end{equation}
and similar equations with $w_1$ and $w_2$ interchanged.
\end{theorem}
\begin{proof}
By \eqref{eqn:Psi-normalization}, it suffices to prove the equations
for $\Psi_{\boldsymbol\gamma}$. Substituting \eqref{eqn:Psi-shift} and its
$-iw_1$ analogue into the definition of
$\mathcal L_{\boldsymbol\gamma,x}^{w_1,w_2}$, and then using
\eqref{eqn:AW-B-identity}, \eqref{eqn:AW-C-identity}, and
Proposition~\ref{prop:psi-algebraic-difference}, gives
\[
\mathcal L_{\boldsymbol\gamma,x}^{w_1,w_2}
\Psi_{\boldsymbol\gamma}(x,y)
=
\left(
-\frac{i(q-q^{-1})}{4}\mu_y^\sigma
+
\cosh\frac{\pi}{w_2}
\left(\rho+\frac{iw_1}{2}-u\pm\sigma\right)
\right)
\Psi_{\boldsymbol\gamma}(x,y).
\]
Using $\mu_y^\sigma=i(\nu_y-\nu_\sigma)/(q-q^{-1})$ and
$i\widehat\gamma_0=\rho-u$, the coefficient is
$v(w_1,w_2,\boldsymbol\gamma;y)$. The duality of $R$ yields the
equation in $y$ with step $iw_1$, and the symmetry
$w_1\leftrightarrow w_2$ yields the two equations with step $iw_2$; see
\cite{R1,R2}.
\end{proof}
For later use, define the conjugated difference operator
\begin{equation*}
	L_{\boldsymbol\gamma,x}^{w_1,w_2}f(x)
	:=
	\frac{
		\mathcal L_{\boldsymbol\gamma,x}^{w_1,w_2}
		\left[G(\pm(\,\cdot\,)-i\gamma_0)f\right](x)}
	{G(\pm x-i\gamma_0)}.
\end{equation*}
Then $\psi_{\boldsymbol\gamma}$ satisfies
\begin{equation}
\begin{split}
	L_{\boldsymbol\gamma,x}^{w_1,w_2}
	\psi_{\boldsymbol\gamma}(x,y)
	&=
	v(w_1,w_2,\boldsymbol\gamma;y)
	\psi_{\boldsymbol\gamma}(x,y),
	\\
	L_{\widehat{\boldsymbol\gamma},y}^{w_1,w_2}
	\psi_{\boldsymbol\gamma}(x,y)
	&=
	v(w_1,w_2,\widehat{\boldsymbol\gamma};x)
	\psi_{\boldsymbol\gamma}(x,y),
\end{split}
\label{eqn:psi-eigenvalue-L}
\end{equation}
and similar equations with $w_1$ and $w_2$ interchanged.

\subsection{Multivariate difference operators}
\label{subsec:multivariate-difference-equations}

We derive a system of difference equations for the multivariate overlap
function \eqref{eqn:multivariate-psi}. The construction is the hyperbolic counterpart of the commuting
multivariate Askey--Wilson operators in \cite{G1,G2}. Because the parameters in
neighboring factors of \eqref{eqn:multivariate-psi} are coupled, the operator
$\mathcal D_j^{w_1,w_2}$ associated with the last $j$ variables, defined below,
contains up to $3^j$ shifts.

In this subsection we assume that
$\boldsymbol u=(u,\ldots,u)$, since the coproduct argument requires the same parameter $u$ in all tensor
factors. We retain the conventions $x_0=\rho$ and
$y_{N+1}=\sigma$ from
\eqref{eqn:multivariate-boundary-parameters}.

The coefficients of the multivariate operators are obtained by repeatedly
applying the one-variable three-term relation from
Proposition~\ref{prop:H-full-step}. Since a shift of $x_k$ changes both
the spectral parameter of the $k$th factor and the boundary parameter of the
$(k+1)$st factor, the coefficient at level $k+1$ depends on the shift chosen
at level $k$. We introduce the notation in the order in which this dependence
arises.

For $k\in\{1,\ldots,N\}$, let $D_k$ denote the backward shift in the
$k$th variable,
\begin{equation}
	(D_k f)(\boldsymbol x)
	:=
	f(x_1,\ldots,x_{k-1},x_k-iw_1,
	x_{k+1},\ldots,x_N).
\end{equation}
For $\boldsymbol\xi=(\xi_1,\ldots,\xi_N)\in\{-1,0,1\}^N$, set
\begin{equation}
	D^{\boldsymbol\xi}
	:=D_1^{\xi_1}\cdots D_N^{\xi_N},
	\qquad
	(D^{\boldsymbol\xi}f)(\boldsymbol x)
	=f(\boldsymbol x-iw_1\boldsymbol\xi).
	\label{eqn:multivariate-combined-shift}
\end{equation}
Thus $\xi_k=-1,0,1$ corresponds, respectively, to shifting $x_k$ by
$+iw_1,0,-iw_1$. The $j$th operator acts only on the final $j$ variables, so
we use the index set
\begin{equation}
	\Xi^{(j)}
	:=\{0\}^{N-j}\times\{-1,0,1\}^j,
	\qquad j=1,\ldots,N.
\end{equation}

Recall from \eqref{eqn:multivariate-A-coefficient} the half-step coefficient
\begin{equation*}
	A_{x,r}^{\lambda}
	:=
	i\frac{
		\cosh\frac{\pi}{w_2}
		\left(x+r+\frac{iw_1}{2}-\lambda\right)
	}{
		\sinh\frac{2\pi x}{w_2}
	},
\end{equation*}
and recall from Proposition \ref{prop:H-full-step} that 
\[
\pi_\lambda(K^{-1})H_{x,r}^{\lambda,u}
=
B_{x,r}^{\lambda}H_{x+iw_1,r}^{\lambda,u}
+C_{x,r}^{\lambda}H_{x,r}^{\lambda,u}
+B_{-x,r}^{\lambda}H_{x-iw_1,r}^{\lambda,u},
\]
with full-step coefficients $B_{x,r}^{\lambda}$ and $C_{x,r}^{\lambda}$ given in terms of $A_{x,r}^{\lambda}$.
These are the coefficients used when the boundary parameter $r$ is unchanged.

We also need the relation in which $r$ is shifted together with $x$. If both are shifted by $+iw_1$, in the same way as in Proposition \ref{prop:H-full-step}, two half-step calculations give
\[
\pi_\lambda(K^{-1})H_{x,r}^{\lambda,u}
=
\widehat B_{x,r}^{\lambda}H_{x+iw_1,r+iw_1}^{\lambda,u}
+\widehat C_{x,r}^{\lambda}H_{x,r+iw_1}^{\lambda,u}
+\widehat B_{-x,r}^{\lambda}H_{x-iw_1,r+iw_1}^{\lambda,u},
\]
with
\begin{gather*}
	\widehat B_{x,r}^{\lambda}
	:=
	A_{x,r}^{\lambda}
	A_{x+iw_1/2,r+iw_1/2}^{\lambda},\\
	\begin{split}
		\widehat C_{x,r}^{\lambda}
		:={}&
		A_{x,r}^{\lambda}
		A_{-x-iw_1/2,r+iw_1/2}^{\lambda}
		+
		A_{-x,r}^{\lambda}
		A_{x-iw_1/2,r+iw_1/2}^{\lambda}.
	\end{split}
\end{gather*}
To treat the three possible shifts of the preceding variable uniformly, for
$\epsilon\in\{-1,0,1\}$ define
\begin{equation}
	B_{x,r}^{\lambda;\epsilon}
	:=
	\begin{cases}
		\widehat B_{x,r}^{\lambda},&\epsilon=-1,\\
		B_{x,r}^{\lambda},&\epsilon=0,\\
		\widehat B_{x,-r}^{\lambda},&\epsilon=1,
	\end{cases}
	\qquad
	C_{x,r}^{\lambda;\epsilon}
	:=
	\begin{cases}
		\widehat C_{x,r}^{\lambda},&\epsilon=-1,\\
		C_{x,r}^{\lambda},&\epsilon=0,\\
		\widehat C_{x,-r}^{\lambda},&\epsilon=1.
	\end{cases}
\end{equation}
The sign convention follows from
\eqref{eqn:multivariate-combined-shift}: $\epsilon=-1$ shifts the preceding
boundary parameter by $+iw_1$, whereas $\epsilon=1$ shifts it by $-iw_1$.

Set $\xi_0:=0$. For a shift vector
$\boldsymbol\xi\in\{-1,0,1\}^N$, the local coefficient at level $k$ is
\begin{equation}
	\mathfrak e_{\boldsymbol\xi,\lambda_k}^{(k)}(\boldsymbol x)
	:=
	\begin{cases}
		B_{x_k,x_{k-1}}^{\lambda_k;\xi_{k-1}},
		&\xi_k=-1,\\[1mm]
		C_{x_k,x_{k-1}}^{\lambda_k;\xi_{k-1}},
		&\xi_k=0,\\[1mm]
		B_{-x_k,x_{k-1}}^{\lambda_k;\xi_{k-1}},
		&\xi_k=1.
	\end{cases}
	\label{eqn:multivariate-local-coefficient}
\end{equation}
The dependence on $\xi_k$ records the shift of $x_k$, while the dependence
on $\xi_{k-1}$ records the induced shift of the boundary parameter
$x_{k-1}$. The coefficient of the complete shift
$D^{\boldsymbol\xi}$ in the last $j$ factors is therefore the product
\begin{equation}
	\boldsymbol{\mathfrak e}_{\boldsymbol\xi,\boldsymbol\lambda}^{(j)}
	(\boldsymbol x)
	:=
	\prod_{k=N-j+1}^{N}
	\mathfrak e_{\boldsymbol\xi,\lambda_k}^{(k)}(\boldsymbol x),
	\qquad \boldsymbol\xi\in\Xi^{(j)}.
	\label{eqn:multivariate-product-coefficient}
\end{equation}
This product is precisely the coefficient obtained by expanding the action of
$(K^{-1})^{\otimes j}$ successively from the factor $N-j+1$ to the factor
$N$.\\

It remains to include the coefficient produced by the relation between
$X_{0,\sigma}$ and $Y_{u,x_{N-j}}$. For this we use the coefficients $a_s$, $b_s^{\rho,\sigma}$ and $\mathfrak b_{u,\epsilon}^{s,r}$ from \eqref{eqn:a-s}, \eqref{eqn:b-s-rho-sigma} and \eqref{eqn:multivariate-b-coefficient}.
Combining the product coefficient
\eqref{eqn:multivariate-product-coefficient} with the terminal coefficient
\eqref{eqn:multivariate-b-coefficient}, we define
\begin{equation}
	\mathcal C_{\boldsymbol\xi}^{(j)}(\boldsymbol x)
	:=
	\mathfrak b_{u,\xi_N}^{x_N,\sigma}
	\boldsymbol{\mathfrak e}_{\boldsymbol\xi,\boldsymbol\lambda}^{(j)}
	(\boldsymbol x)
	-
	\delta_{\boldsymbol\xi,\boldsymbol0}
	b_u^{x_{N-j},\sigma}.
\end{equation}
The Kronecker-delta term comes from the $-1$ in
$\boldsymbol K^{-1,(j)}-1$ and therefore contributes only to the unshifted
term. The $j$th multivariate difference operator is
\begin{equation}
	\mathcal D_j^{w_1,w_2}
	:=
	\sum_{\boldsymbol\xi\in\Xi^{(j)}}
	\mathcal C_{\boldsymbol\xi}^{(j)}(\boldsymbol x)
	D^{\boldsymbol\xi}.
	\label{eqn:multivariate-difference-operator}
\end{equation}
The operator $\mathcal D_j^{w_1,w_2}$ shifts only the final $j$ components of $\boldsymbol x$ and contains at most
$3^j$ shifts. Theorem \ref{thm:multivariate-Ruijsenaars-difference-equations} will establish the corresponding eigenvalue equation for
the multivariate overlap function.

In the notation of \eqref{eqn:K-Y-tau-operators}, let
\[
\mathsf K^{(j)}:=\sum_{\boldsymbol\xi\in\Xi^{(j)}}\boldsymbol{\mathfrak e}_{\boldsymbol\xi,\boldsymbol\lambda}^{(j)}(\boldsymbol x)D^{\boldsymbol\xi},
\]
and let $\mathsf Y^{(j)}$ denote multiplication by $\mu_{x_N}^{x_{N-j}}$. Since $D^{\boldsymbol\xi}$ does not shift $x_{N-j}$ for $\boldsymbol\xi\in\Xi^{(j)}$, identity \eqref{eqn:mu-b-identity}, with $(\rho,x,\epsilon)$ replaced by $(x_{N-j},x_N,\xi_N)$, shows that
\begin{equation}
	\mathcal D_j^{w_1,w_2}
	=
	a_u\mathsf Y^{(j)}\mathsf K^{(j)}
	+a_{-u}\mathsf K^{(j)}\mathsf Y^{(j)}
	+b_u^{x_{N-j},\sigma}\bigl(\mathsf K^{(j)}-1\bigr),
	\label{eqn:multivariate-D-operator-form}
\end{equation}
which is the multivariate analogue of \eqref{eqn:D-operator}. For $N=j=1$ we recover $\mathcal D_1^{w_1,w_2}=\mathsf D_{\rho,\sigma}^{\lambda_1,u}$.

For completeness, let us give the explicit formulas for the relevant coefficients in terms of hyperbolic and trigonometric functions (the coefficients $\mathfrak b_{u,\epsilon}^{s,r}$ are given in \eqref{eqn:b-closed-form}):
\begin{equation}
	\begin{split}
		B_{x,r}^{\lambda}
		&=
		\frac{
			\cosh\frac{\pi}{w_2}
			\left(x+r-\lambda+\frac{iw_1}{2}\right)
			\cosh\frac{\pi}{w_2}
			\left(x-r-\lambda+\frac{iw_1}{2}\right)
		}{
			\sinh\frac{2\pi x}{w_2}
			\sinh\frac{2\pi(x+iw)}{w_2}
		},\\
		\widehat B_{x,r}^{\lambda}
		&=-
		\frac{
			\cosh\frac{\pi}{w_2}
			\left(x+r+\frac{iw_1}{2}-\lambda\right)
			\cosh\frac{\pi}{w_2}
			\left(x+r+\frac{3iw_1}{2}-\lambda\right)
		}{
			\sinh\frac{2\pi x}{w_2}
			\sinh\frac{2\pi(x+iw_1/2)}{w_2}
		},\\
		C_{x,r}^{\lambda}
		&=-
		\frac{
			\cosh\frac{\pi}{w_2}
			\left(x+r+\lambda+\frac{iw_1}{2}\right)
			\cosh\frac{\pi}{w_2}
			\left(x+r-\lambda+\frac{iw_1}{2}\right)
		}{
			\sinh\frac{2\pi x}{w_2}
			\sinh\frac{2\pi(x+iw)}{w_2}
		}
		\\
		&\qquad -
		\frac{
			\cosh\frac{\pi}{w_2}
			\left(x-r+\lambda-\frac{iw_1}{2}\right)
			\cosh\frac{\pi}{w_2}
			\left(x-r-\lambda-\frac{iw_1}{2}\right)
		}{
			\sinh\frac{2\pi x}{w_2}
			\sinh\frac{2\pi(x-iw)}{w_2}
		},\\
		\widehat C_{x,r}^{\lambda} 
		&= 
		\frac{ 2\cos\left(\frac{\pi w_1}{w_2}\right) \cosh\frac{\pi}{w_2} \left(x+r+\frac{iw_1}{2}-\lambda\right) \cosh\frac{\pi}{w_2} \left(x-r-\frac{iw_1}{2}+\lambda\right) }{ \sinh\frac{2\pi}{w_2} \left(x+\frac{iw_1}{2}\right) \sinh\frac{2\pi}{w_2} \left(x-\frac{iw_1}{2}\right) }.
	\end{split} 
\end{equation}

The following lemma shows that the product
$\boldsymbol{\mathfrak e}_{\boldsymbol\xi,\boldsymbol\lambda}^{(j)}$ is exactly the
coefficient in the expansion of the difference operator corresponding to iterated coproducts of $K^{-1}$.
\begin{lemma}
	\label{lem:multivariate-K-inverse-action}
	Let $j\in\{1,\ldots,N\}$, and set
	\begin{equation}
		\boldsymbol K^{-1,(j)}
		:=
		1^{\otimes(N-j)}\otimes(K^{-1})^{\otimes j}.
	\end{equation}
	Away from the polar hyperplanes of the coefficients above,
	\begin{equation}
		\pi_{\boldsymbol\lambda}
		\bigl(\boldsymbol K^{-1,(j)}\bigr)
		\mathbf H_{\boldsymbol x,\rho}^{\boldsymbol\lambda,u}
		=
		\sum_{\boldsymbol\xi\in\Xi^{(j)}}
		\boldsymbol{\mathfrak e}_{\boldsymbol\xi,\boldsymbol\lambda}^{(j)}(\boldsymbol x)
		D^{\boldsymbol\xi}
		\mathbf H_{\boldsymbol x,\rho}^{\boldsymbol\lambda,u}.
		\label{eqn:multivariate-K-inverse-expansion}
	\end{equation}
\end{lemma}

\begin{proof}
	Apply the one-variable three-term relation first to the factor with index
	$N-j+1$. This produces the three possibilities
	$\xi_{N-j+1}=-1,0,1$. A nonzero shift of $x_{N-j+1}$ also changes the
	second parameter of the next factor. Consequently, at the next step one uses
	$\widehat B_{x,r}$ and $\widehat C_{x,r}$ when
	$\xi_{N-j+1}=-1$, $\widehat B_{x,-r}$ and $\widehat C_{x,-r}$
	when $\xi_{N-j+1}=1$, and the unshifted coefficients $B$ and $C$ when
	$\xi_{N-j+1}=0$. This is precisely the dependence on
	$\xi_{k-1}$ in \eqref{eqn:multivariate-local-coefficient}.
	
	Continuing through the factors $k=N-j+2,\ldots,N$ produces one term for each
	$\boldsymbol\xi\in\Xi^{(j)}$. The product of the local coefficients is
	\eqref{eqn:multivariate-product-coefficient}, and the combined shift is
	$D^{\boldsymbol\xi}$. This proves
	\eqref{eqn:multivariate-K-inverse-expansion}.
\end{proof}

We can now state the multivariate difference equations. 
\begin{theorem}
	\label{thm:multivariate-Ruijsenaars-difference-equations}
Assume that $u_1=\cdots=u_N=u$, and let the operators
$\mathcal D_j^{w_1,w_2}$ be defined by
\eqref{eqn:multivariate-difference-operator}.
The multivariate overlap function is a joint eigenfunction
of the family $\mathcal D_1^{w_1,w_2},\ldots,\mathcal D_N^{w_1,w_2}$, satisfying
\begin{equation}
	\mathcal D_j^{w_1,w_2}
	\boldsymbol\psi_{\rho,\sigma}^{\boldsymbol\lambda,u}
	(\boldsymbol x,\boldsymbol y)
	=
	\mu_{y_{N-j+1}}^\sigma
	\boldsymbol\psi_{\rho,\sigma}^{\boldsymbol\lambda,u}
	(\boldsymbol x,\boldsymbol y),
	\qquad j=1,\ldots,N.
	\label{eqn:multivariate-Ruijsenaars-difference-equation}
\end{equation}
These equations hold as identities of meromorphic functions.
\end{theorem}

\begin{proof}
	We first assume $w_2>3w_1>0$, $\rho,\sigma,x_k,y_k\in\R$ and $\lambda_k,u\in i(-\delta,\delta)$; the general case then follows by meromorphic continuation, see Proposition~\ref{prop:difference-equations-complex-periods}. For these parameters we apply Lemma~\ref{lem:straight-contour} to the $k$th factor of \eqref{eqn:multivariate-psi}, with $(x_{k-1},y_{k+1})$ in the role of $(\rho,\sigma)$. This gives horizontal contours $\mathcal C_k$. The shifts occurring below change $x_k$ and $x_{k-1}$ by at most $iw_1$, so by Lemma~\ref{lem:straight-contour}(ii) the same contours $\mathcal C_k$ compute all shifted multivariate overlaps.
	Recall from \eqref{eqn:multivariate-F-eigenvalue-positive} that the multivariate
	$F$-function is an eigenfunction of $\pi_{\overline{\boldsymbol\lambda}}\bigl((\mathbf X_{\sigma}^{(j)})^*\bigr)$, with eigenvalue
	$\overline{\mu_{y_{N-j+1}}^{\sigma}}$. Applying Lemma~\ref{lem:sesquilinear-star-compatibility} in each tensor factor, which is allowed by Lemma~\ref{lem:straight-contour}(i), therefore gives
	\begin{equation}
		\mu_{y_{N-j+1}}^{\sigma}
		\boldsymbol\psi_{\rho,\sigma}^{\boldsymbol\lambda,u}
		(\boldsymbol x,\boldsymbol y)
		=
		\left\langle
		\pi_{\boldsymbol\lambda}(\mathbf X_{\sigma}^{(j)})
		\mathbf H_{\boldsymbol x,\rho}^{\boldsymbol\lambda,u},
		\mathbf F_{\boldsymbol y,\sigma}^{\boldsymbol\lambda,0}
		\right\rangle.
		\label{eqn:move-X-through-multivariate-pairing}
	\end{equation}
	Apply the coproduct to the relation \eqref{eqn:relation-XY} with $\rho=x_{N-j}$. Since the coproduct is an algebra homomorphism, we
	obtain
	\begin{equation}
		\mathbf X_{\sigma}^{(j)}
		=
		a_u\boldsymbol K^{-1,(j)}
		\mathbf Y_{u,x_{N-j}}^{(j)}
		+
		a_{-u}\mathbf Y_{u,x_{N-j}}^{(j)}
		\boldsymbol K^{-1,(j)}
		+
		b_u^{x_{N-j},\sigma}
		\bigl(\boldsymbol K^{-1,(j)}-1\bigr).
		\label{eqn:multivariate-X-Y-relation}
	\end{equation}
	By Proposition~\ref{prop:multivariate-H-eigenfunctions} and Lemma~\ref{lem:multivariate-K-inverse-action},
	\[
	\pi_{\boldsymbol\lambda}\bigl(\mathbf Y_{u,x_{N-j}}^{(j)}\bigr)\mathbf H_{\boldsymbol x,\rho}^{\boldsymbol\lambda,u}=\mathsf Y^{(j)}\mathbf H_{\boldsymbol x,\rho}^{\boldsymbol\lambda,u},
	\qquad
	\pi_{\boldsymbol\lambda}\bigl(\boldsymbol K^{-1,(j)}\bigr)\mathbf H_{\boldsymbol x,\rho}^{\boldsymbol\lambda,u}=\mathsf K^{(j)}\mathbf H_{\boldsymbol x,\rho}^{\boldsymbol\lambda,u}.
	\]
	The operators $\mathsf Y^{(j)}$ and $\mathsf K^{(j)}$ act in $\boldsymbol x$ and do not shift $x_{N-j}$, so they commute with the operators $\pi_{\boldsymbol\lambda}(m)$ occurring here. As in the one-variable case the order of the factors is reversed, and \eqref{eqn:multivariate-X-Y-relation} and \eqref{eqn:multivariate-D-operator-form} give
	\[
	\pi_{\boldsymbol\lambda}(\mathbf X_{\sigma}^{(j)})\mathbf H_{\boldsymbol x,\rho}^{\boldsymbol\lambda,u}
	=\mathcal D_j^{w_1,w_2}\mathbf H_{\boldsymbol x,\rho}^{\boldsymbol\lambda,u}.
	\]
	The right-hand side is a linear combination of the functions $D^{\boldsymbol\xi}\mathbf H_{\boldsymbol x,\rho}^{\boldsymbol\lambda,u}$, $\boldsymbol\xi\in\Xi^{(j)}$, with coefficients independent of $\boldsymbol z$. Substituting into \eqref{eqn:move-X-through-multivariate-pairing}, integrating term by term, and using the contours $\mathcal C_k$ chosen above for the shifted overlaps proves
	\eqref{eqn:multivariate-Ruijsenaars-difference-equation}.
\end{proof}

\begin{remark}
\label{rem:commutativity-formal}
Formally, the operators $\mathcal D_j^{w_1,w_2}$ commute because the algebra elements $\mathbf X_\sigma^{(j)}$ commute by \eqref{eqn:commuting-X-family}. We do not need this property in the sequel. To make it rigorous one may argue as follows: the commutator $[\mathcal D_j^{w_1,w_2},\mathcal D_k^{w_1,w_2}]$ is a difference operator with shifts in $\{-2,\ldots,2\}^N$ that annihilates $\boldsymbol\psi(\cdot,\boldsymbol y)$ for all $\boldsymbol y$, and the shifted functions $\boldsymbol\psi(\boldsymbol x-iw_1\boldsymbol\xi,\cdot)$ are linearly independent for generic $\boldsymbol x$, since by the reversal symmetry \eqref{eqn:multivariate-psi-reversal} and the eigenvalue equations they are joint eigenfunctions in $\boldsymbol y$ with pairwise distinct eigenvalues.
\end{remark}

Interchanging $w_1$ and $w_2$ gives a second system of $N$ equations with
shifts by $iw_2$. Applying the reversal symmetry
\eqref{eqn:multivariate-psi-reversal} gives two dual systems acting on
$\boldsymbol y$. Hence the multivariate Ruijsenaars function satisfies
four systems of $N$ difference equations, generalizing the four
one-variable Askey--Wilson equations.

Instead of using the reversal symmetry, the dual difference equations acting
in the $\boldsymbol y$-variables can also be obtained
directly from the modular-double representation. The automorphism
$\theta$, which satisfies $\theta(Y_{u,\rho})=X_{u,\rho}$, reverses the coproduct. Hence, after composing with the reversal of
the tensor factors, the nested coproduct family associated with
$Y_{u,\rho}$ gives a nested commuting family associated with
$X_{u,\rho}$. The multivariate $F$-function is a joint eigenfunction
of the appropriate adjoints of this family. Transferring their action
through the multivariate overlap pairing and applying the three-term
relation for the $F$-functions therefore yields difference operators
acting in $\boldsymbol y$. 	

\subsection{Extension to complex \texorpdfstring{$w_1,w_2$}{w1, w2}}
\label{subsec:difference-equations-complex-periods}

The difference equations remain valid with weaker conditions on the parameters $w_1$ and $w_2$. Recall the domain $\mathcal P$ from Section~\ref{subsec:notation}.
On this domain $\arg(w_1w_2)\in(-\pi,\pi)$, and we use the principal branch of $\sqrt{w_1w_2}$; it is positive for positive $w_1,w_2$ and equals $|w_1|$ for a complex conjugate pair. Define the continued overlap kernel by
\eqref{eqn:Rren-psi-relation}, and the multivariate kernel by the product
\eqref{eqn:multivariate-psi}.

\begin{proposition}
\label{prop:difference-equations-complex-periods}
The difference equations of Proposition~\ref{prop:psi-algebraic-difference} and Theorem~\ref{thm:R-askey-wilson}, and the four multivariate systems
of Theorem~\ref{thm:multivariate-Ruijsenaars-difference-equations} hold
for $(w_1,w_2)\in\mathcal P$ satisfying $q^4\neq1, \widetilde q^{\,4}\neq1$. They are identities of meromorphic functions in the remaining parameters
and variables. 
\end{proposition}

\begin{proof}
The restrictions on $q$ and $\widetilde q$ ensure that $q-q^{-1}$, $q^2-q^{-2}$, and their tilded counterparts are nonzero.
The joint meromorphic continuation of Ruijsenaars' function and the
hyperbolic gamma function in $w_1,w_2$, the parameters, and the variables
\cite[Theorem~2.2]{R1} makes each continued overlap kernel meromorphic
on $\mathcal P$ and in its remaining arguments. The same holds for its
finite products, all coefficients, and all shifted kernels occurring
in the equations. Consider, for example, the $j$th equation of Theorem~\ref{thm:multivariate-Ruijsenaars-difference-equations}, and let
\[
F_j:=
\mu_{y_{N-j+1}}^\sigma
\boldsymbol\psi_{\rho,\sigma}^{\boldsymbol\lambda,u}
(\boldsymbol x,\boldsymbol y)
-\sum_{\boldsymbol\xi\in\Xi^{(j)}}
\mathcal C_{\boldsymbol\xi}^{(j)}(\boldsymbol x)
\boldsymbol\psi_{\rho,\sigma}^{\boldsymbol\lambda,u}
(\boldsymbol x-iw_1\boldsymbol\xi,\boldsymbol y),
\]
regarded as a function of $(w_1,w_2)$ and all remaining parameters and variables. It is meromorphic on a connected domain of the form $\mathcal P\times\C^m$. By the proofs above, $F_j$ vanishes for $w_2>3w_1>0$ with $w_2/w_1\notin\frac12\Z$, $\rho,\sigma,x_k,y_k\in\R$ and $\lambda_k,u\in i(-\delta,\delta)$ with $\delta=\frac{w_2-3w_1}{2}$. This set contains a product of nonempty open intervals, one in each variable, so $F_j$ vanishes identically by the identity theorem applied in one variable at a time. The same argument applies to the one-variable equations and to the other three multivariate systems; for the systems with step $iw_2$ one interchanges the roles of $w_1$ and $w_2$, and the systems in $\boldsymbol y$ follow from the reversal symmetry \eqref{eqn:multivariate-psi-reversal}.
\end{proof}

The coefficients of the equations \eqref{eqn:psi-eigenvalue-L} in Askey--Wilson form depend meromorphically on $(w_1,w_2)\in\mathcal P$, but, unlike those of Proposition~\ref{prop:psi-algebraic-difference}, they contain no factors $(q-q^{-1})^{-1}$ or $(q^2-q^{-2})^{-1}$ and are not singular along $q^4=1$ or $\widetilde q^{\,4}=1$. Hence the same argument shows that these equations hold for all $(w_1,w_2)\in\mathcal P$.

In Section~\ref{sec:ruijsenaars-transforms} we again restrict to positive $w_1,w_2$ or complex conjugate pairs, which is needed for the positivity of the weight functions and for the adjoint arguments; the condition $q^4\neq1,\widetilde q^{\,4}\neq 1$ is not needed there.

\section{The Ruijsenaars transform}
\label{sec:ruijsenaars-transforms}

We now establish the transform associated with the overlap kernel
$\psi_{\boldsymbol\gamma}$. Ruijsenaars proved the corresponding result for
positive $w_1,w_2$ and real $\boldsymbol\gamma$ in \cite{R3}. Here we use the
parameter regime needed for the multivariate construction: $\rho,\sigma$ are
real, $u,\lambda$ are purely imaginary, and $w_1,w_2$ may also form a
nonreal complex conjugate pair. Ruijsenaars' proof in \cite{R3} uses time-dependent scattering theory. Here we follow instead the Green identity and asymptotic method used for the Askey--Wilson function transform in \cite{KS}. The results of this section hold for all positive $w_1,w_2$ and for all complex conjugate pairs; no condition on $q$ is needed. Indeed, the eigenvalue equation \eqref{eqn:psi-eigenvalue-L}, derived above under the assumption $q^4\neq1, \widetilde q^{\,4}\neq 1$, has coefficients that are not singular along $q^4=1$ or $\widetilde q^{\,4}=1$, and it holds for all such $w_1,w_2$ by Section~\ref{subsec:difference-equations-complex-periods}; see also \cite{R1}.

Set
\[
w_0:=
\begin{cases}
	\min\{w_1,w_2\},&w_1,w_2>0,\\
	w,&\overline{w_1}=w_2.
\end{cases}
\]
In the positive case we may interchange $w_1$ and $w_2$ and assume that
$0<w_1=w_0<w_2$. We refer to this regime as the \emph{positive case}; the
regime $\overline{w_1}=w_2$, with $\operatorname{Re}w_1>0$ and
$\operatorname{Im}w_1\neq0$, is called the \emph{conjugate case}.

\subsection{The weighted transform}

Define
\begin{equation}
	c'(\boldsymbol\gamma;x)
	:=
	\frac{
		G(x+i\gamma_0)\prod_{j=1}^3G(x-i\gamma_j)
	}{
		G(2x+iw)
	}
	\label{eqn:c-prime}
\end{equation}
and
\begin{equation}
	\begin{split}
	W_{\boldsymbol\gamma}(x)
	&:=
	\frac{1}{w_1w_2}
	\frac{1}{c'(\boldsymbol\gamma;x)c'(\boldsymbol\gamma;-x)}.
	\end{split} \label{eqn:transform-weight}
\end{equation}
In terms of $(\lambda,u,\rho,\sigma)$, we have
\[
W_{\boldsymbol\gamma}(x)
= \frac{1}{w_1w_2}
\frac{ G(iw\pm 2x)}{G(\lambda\pm\rho \pm x)
	G(u\pm\sigma \pm x)}.
\]
If $u,\lambda\in i\R$ with $|\operatorname{Im}u|,|\operatorname{Im}\lambda|<w$ and $\rho,\sigma\in\R$, then
$W_{\boldsymbol\gamma}(x)>0$ for $x>0$. Indeed, all arguments of $G$ in the denominator then lie in the strip $|\operatorname{Im}z|<w$, where $G$ is analytic and zero-free. By \eqref{eqn:G-symmetries} (note that $G(w_1,w_2;\cdot)=G(\overline{w_1},\overline{w_2};\cdot)$ in the conjugate case, since $G$ is symmetric in $w_1,w_2$), $\overline{G(z)}=G(-\overline z)$, so that for real $a$ and $x$ and $\lambda\in i\R$ we have $G(\lambda+a+x)G(\lambda-a-x)=|G(\lambda+a+x)|^2>0$, and similarly for $u$. Hence
$G(\lambda\pm\rho\pm x)=|G(\lambda+\rho+x)|^2|G(\lambda-\rho+x)|^2>0$, and similarly $G(u\pm\sigma\pm x)>0$. Finally $G(iw\pm2x)=|G(iw+2x)|^2>0$ for $x>0$, because the only zero of $G$ on the line $\operatorname{Im}z=w$ is $z=iw$ and the poles of $G$ lie in $\operatorname{Im}z\leq-w$. We write
\[
\mathcal H_{\boldsymbol\gamma}
:=L^2\bigl((0,\infty),W_{\boldsymbol\gamma}(x)\dd x\bigr)
\]
and
\[
\langle f,g\rangle_{\boldsymbol\gamma}
:=
\int_0^\infty f(x)\overline{g(x)}W_{\boldsymbol\gamma}(x)\dd x.
\]
For $f\in C_{\mathrm c}^\infty(0,\infty)$, define
\begin{equation}
	(\mathcal R_{\boldsymbol\gamma}f)(y)
	:=
	\int_0^\infty
	f(x)\psi_{\boldsymbol\gamma}(x,y)
	W_{\boldsymbol\gamma}(x)\dd x.
	\label{eqn:weighted-R-transform}
\end{equation}
We show that under suitable parameter conditions, the transform $\mathcal R_{\boldsymbol \gamma}$ extends to a unitary map
$\mathcal H_{\boldsymbol\gamma}\to\mathcal H_{\widehat{\boldsymbol\gamma}}$.

\subsection{Asymptotics of the kernel}

The transform proof requires an asymptotic expansion of
$\psi_{\boldsymbol\gamma}(x,y)$ as
$\operatorname{Re}x\to\infty$. For real parameters this is equivalent to
the asymptotic expansion of Ruijsenaars' normalized function derived in
\cite{R2}. The following version also allows complex parameters and conjugate
$w_1,w_2$.

Set
\begin{equation}
	\begin{split}
		\psi_{\boldsymbol\gamma}^{\mathrm{as}}(x,y)
		:={}&
		\sqrt{w_1w_2}\,
		e^{-\alpha(w+iu+i\lambda)x}
		\left(
		c'(\widehat{\boldsymbol\gamma};y)e^{i\alpha y x}
		+
		c'(\widehat{\boldsymbol\gamma};-y)e^{-i\alpha y x}
		\right),
	\end{split}
\end{equation}
where $c'$ is given by \eqref{eqn:c-prime}.

\begin{theorem}
	\label{thm:psi-asymptotics-article}
	Let $\lambda,u\in i\R$ with
	\[
	\operatorname{Im}\lambda<w,
	\qquad
	\operatorname{Im}u<w,
	\qquad
	\operatorname{Im}(\lambda+u)>w_0-w,
	\]
	let $\rho,\sigma\in\R$, and let $0<r_1<r_2$. Then there exist $C,R,\varepsilon>0$ such that
	\begin{equation}
		\left|
		\psi_{\boldsymbol\gamma}(x,y)
		-
		\psi_{\boldsymbol\gamma}^{\mathrm{as}}(x,y)
		\right|
		\leq
		C e^{-\alpha(w-\operatorname{Im}(\lambda+u)+\varepsilon)\operatorname{Re}x}
		\label{eqn:psi-asymptotic-estimate}
	\end{equation}
	for $\operatorname{Re}x>R$, $|\operatorname{Im}x|\leq w_0$ and $y\in[r_1,r_2]$.
\end{theorem}
Since $|e^{-\alpha(w+iu+i\lambda)x}|=e^{-\alpha(w-\operatorname{Im}(\lambda+u))\operatorname{Re}x}$, the leading term satisfies $\psi_{\boldsymbol\gamma}^{\mathrm{as}}(x,y)=O\bigl(e^{-\alpha(w-\operatorname{Im}(\lambda+u))\operatorname{Re}x}\bigr)$, so \eqref{eqn:psi-asymptotic-estimate} says that the error is smaller than the leading term by a factor $e^{-\alpha\varepsilon\operatorname{Re}x}$.

\begin{proof}
	After shifting the integration variable, we have
	\[
	\psi_{\boldsymbol\gamma}(x,y)
	=
	\int_{\mathcal C}
	\frac{G(z+\lambda\pm x)}
	{G(z+iw\pm\rho)}
	\frac{G(z-u\pm y)}
	{G(z+iw+\lambda-u\pm\sigma)}
	\,\dd z.
	\]
	Fix $\eta\in(0,w_0)$ and choose a contour $\mathcal C_\eta$ below $\mathcal C$ such that the
	deformation from $\mathcal C$ to $\mathcal C_\eta$ crosses precisely
	the two poles
	\[
	z_\pm=u\pm y-iw
	\]
	and no other poles. The parameter assumptions ensure that this downward
	deformation, and the resulting contour, may be chosen uniformly for
	$y\in[r_1,r_2]$.
	
	Using
	\[
	\operatorname{Res}_{z=-iw}G(z)
	=
	\frac{i\sqrt{w_1w_2}}{2\pi},
	\]
	the residue theorem gives
	\begin{equation}
		\begin{split}
			\psi_{\boldsymbol\gamma}(x,y)
			={}&
			\sqrt{w_1w_2}\,
			c'(\widehat{\boldsymbol\gamma};y)
			G(-y-iw+u+\lambda\pm x)
			\\
			&+
			\sqrt{w_1w_2}\,
			c'(\widehat{\boldsymbol\gamma};-y)
			G(y-iw+u+\lambda\pm x)
			\\
			&+
			I_{\mathrm{rem}}(x,y),
		\end{split}
		\label{eqn:psi-residue-decomposition}
	\end{equation}
	where
	\[
	I_{\mathrm{rem}}(x,y)
	:=
	\int_{\mathcal C_\eta}
	\frac{G(z+\lambda\pm x)}
	{G(z+iw\pm\rho)}
	\frac{G(z-u\pm y)}
	{G(z+iw+\lambda-u\pm\sigma)}
	\,\dd z.
	\]	
	For example, the residue at $z=u+y-iw$ comes from the pole of
	$G(z-u-y)$ at $-iw$; simplifying the remaining factors using the
	symmetries of $G$ gives the second residue term. The pole at
	$z=u-y-iw$ similarly gives the first term.
	The uniform asymptotic expansion of the hyperbolic gamma function
	gives, for $y\in[r_1,r_2]$ and
	$|\operatorname{Im}x|\leq w_0$,
	\[
	\begin{split}
	&G(-\epsilon y-iw+u+\lambda+x)\,G(-\epsilon y-iw+u+\lambda-x)\\
	&\qquad=
	e^{-\alpha(w+iu+i\lambda)x}
	e^{i\epsilon\alpha y x}
	\left(1+
	O\left(e^{-\alpha\delta_1\operatorname{Re}x}\right)\right),
	\qquad \epsilon\in\{1,-1\},
	\end{split}
	\]
	for a positive $\delta_1$ that may be chosen uniformly on compact
	subsets of the parameter domain. Hence the residue contribution in
	\eqref{eqn:psi-residue-decomposition} equals
	$\psi_{\boldsymbol\gamma}^{\mathrm{as}}(x,y)$ up to an
	exponentially small error.
	
	It remains to estimate $I_{\mathrm{rem}}$. Divide
	$\mathcal C_\eta$ into its two tails, a compact middle part, and the
	indentations around the moving pole sequences. The uniform
	hyperbolic-gamma estimates from
	Appendix~\ref{sec:hyperbolic-gamma} give
	\[
	|I_{\mathrm{rem}}(x,y)|
	\leq
	C_1e^{-\alpha(2w-w_0)\operatorname{Re}x}
	+
	C_2e^{-\alpha(
		w+\eta-\operatorname{Im}\lambda-\operatorname{Im}u)
		\operatorname{Re}x}.
	\]
	Since $2w-w_0=w-\operatorname{Im}(\lambda+u)+\bigl(w-w_0+\operatorname{Im}(\lambda+u)\bigr)$ and $w+\eta-\operatorname{Im}\lambda-\operatorname{Im}u=w-\operatorname{Im}(\lambda+u)+\eta$, this gives
	\[
	|I_{\mathrm{rem}}(x,y)|
	\leq
	C e^{-\alpha(w-\operatorname{Im}(\lambda+u)+\varepsilon')\operatorname{Re}x},
	\qquad
	\varepsilon':=\min\bigl\{\eta,\ w-w_0+\operatorname{Im}(\lambda+u)\bigr\}>0,
	\]
	where $\varepsilon'>0$ by the assumption $\operatorname{Im}(\lambda+u)>w_0-w$. The error in the residue contribution is of order $e^{-\alpha(w-\operatorname{Im}(\lambda+u)+\delta_1)\operatorname{Re}x}$. Hence \eqref{eqn:psi-asymptotic-estimate} holds with $\varepsilon:=\min\{\varepsilon',\delta_1\}$.
\end{proof}

\subsection{Green's identity and the Casorati determinant}
 
For fixed $y>0$, let
\[
\psi_{\boldsymbol\gamma}^{y}(x):=
\psi_{\boldsymbol\gamma}(x,y).
\]
These functions are even in $x$ and, by the first equation in \eqref{eqn:psi-eigenvalue-L}, satisfy
\begin{equation}
	L_{\boldsymbol\gamma,x}^{w_1,w_2}
	\psi_{\boldsymbol\gamma}^{y}
	=
	v(w_1,w_2,\boldsymbol\gamma;y)
	\psi_{\boldsymbol\gamma}^{y}.
	\label{eqn:psi-eigenvalue-transform-section}
\end{equation}

In the positive case, define the symmetric truncated pairing by
\begin{equation}
	\langle f,g\rangle_{\boldsymbol\gamma,M}
	:=
	\int_{-M}^{M}f(x)\overline g(x)W_{\boldsymbol\gamma}(x)\dd x.
\end{equation}
For even $f,g$, this is twice the integral over $(0,M)$.

In the conjugate case, let $0<\epsilon<|\operatorname{Im}w_1|$ and
$M>2|\operatorname{Im}w_1|$. Let $\mathcal C_{M,\epsilon}$ be obtained from
$[-M,M]$ by replacing the intervals of radius $\epsilon$ around
$\pm\operatorname{Im}w_1$ by semicircles, passing above
$-\operatorname{Im}w_1$ and below $\operatorname{Im}w_1$. We set
\begin{equation}
	\langle f,g\rangle_{\boldsymbol\gamma,M,\epsilon}
	:=
	\int_{\mathcal C_{M,\epsilon}}
	f(x)\overline g(x)W_{\boldsymbol\gamma}(x)\dd x.
\end{equation}
The indentations avoid the poles encountered when the contour is shifted by
$iw_1$.

Define
\begin{equation}
	t_{\boldsymbol\gamma}(x)
	:=
	\frac{
		\cosh\frac{\pi}{w_2}
		\left(x+\frac{iw_1}{2}-i\gamma_0\right)
	}{
		\cosh\frac{\pi}{w_2}
		\left(x+\frac{iw_1}{2}+i\gamma_0\right)
	}.
\end{equation}
The difference equation for $G$ gives
\begin{equation}
	t_{\boldsymbol\gamma}(x)
	=
	G(\pm x+i\gamma_0)
	G\bigl(\pm(x+iw_1)-i\gamma_0\bigr).
\end{equation}
For functions analytic in a strip containing the line segment from
$z-iw_1$ to $z$, define the Casorati determinant
\begin{equation}
	\begin{split}
		[f,g]_{\boldsymbol\gamma}(z)
		:={}&
		\int_{z-iw_1}^{z}
		\Bigl(
		f(x+iw_1)\overline g(x)
		-
		f(x)\overline g(x+iw_1)
		\Bigr)
		\\
		&\qquad\times
		t_{\boldsymbol\gamma}(x)
		\mathcal A_{\boldsymbol\gamma}^{w_1,w_2}(x)
		W_{\boldsymbol\gamma}(x)\dd x.
	\end{split}
\end{equation}
The product
$t_{\boldsymbol\gamma}(x)
\mathcal A_{\boldsymbol\gamma}^{w_1,w_2}(x)
W_{\boldsymbol\gamma}(x)$ is invariant under
$x\mapsto-x-iw_1$. Hence, for even $f$ and $g$,
\begin{equation}
	[f,g]_{\boldsymbol\gamma}(-z)
	=
	-[f,g]_{\boldsymbol\gamma}(z).
	\label{eqn:Casorati-bracket-odd}
\end{equation}

\begin{lemma} 
	Let $f$ and $g$ be even, analytic in a strip containing all contours and
	their $iw_1$-translates used below, and have no poles on or between these
	contours. Assume also sufficient decay for the contour shifts and boundary
	integrals below. In the positive case,
	\begin{equation}
		\begin{split}
			&
			\left\langle
			L_{\boldsymbol\gamma,x}^{w_1,w_2}f,g
			\right\rangle_{\boldsymbol\gamma,M}
			-
			\left\langle
			f,L_{\overline{\boldsymbol\gamma},x}^{w_1,w_2}g
			\right\rangle_{\boldsymbol\gamma,M}
			=
			2[f,g]_{\boldsymbol\gamma}(M).
		\end{split}
		\label{eqn:Green-identity-positive}
	\end{equation}
	In the conjugate case,
	\begin{equation}
		\begin{split}
			&
			\left\langle
			L_{\boldsymbol\gamma,x}^{w_1,w_2}f,g
			\right\rangle_{\boldsymbol\gamma,M,\epsilon}
			-
			\left\langle
			f,L_{\overline{\boldsymbol\gamma},x}^{w_2,w_1}g
			\right\rangle_{\boldsymbol\gamma,M,\epsilon}
			=
			2[f,g]_{\boldsymbol\gamma}(M)
		\end{split}
		\label{eqn:Green-identity-conjugate}
	\end{equation}
	for $M$ sufficiently large and $\epsilon>0$ sufficiently small.
\end{lemma}

\begin{proof}
	We give the positive case. After inserting the definitions of the conjugated
	Askey--Wilson operator and the weight, the left-hand side of
	\eqref{eqn:Green-identity-positive} becomes
	\[
	\left(
	\int_{-M}^{M}
	-
	\int_{-M-iw_1}^{M-iw_1}
	\right)
	\Bigl(
	f(x+iw_1)\overline g(x)
	-
	f(x)\overline g(x+iw_1)
	\Bigr)
	\]
	\[
	{}\times
	t_{\boldsymbol\gamma}(x)
	\mathcal A_{\boldsymbol\gamma}^{w_1,w_2}(x)
	W_{\boldsymbol\gamma}(x)\dd x.
	\]
	Cauchy's theorem reduces this difference to
	$[f,g]_{\boldsymbol\gamma}(M)-[f,g]_{\boldsymbol\gamma}(-M)$, which equals
	$2[f,g]_{\boldsymbol\gamma}(M)$ by
	\eqref{eqn:Casorati-bracket-odd}. In the conjugate case, complex conjugation
	interchanges $w_1$ and $w_2$, and the indented contour ensures that no poles
	are crossed. Deforming the indented contour gives the same boundary
	contribution at $\pm M$. The same boundary calculation proves
	\eqref{eqn:Green-identity-conjugate}.
\end{proof}

Before applying the lemma, we record how the eigenvalue equation interacts with the adjoint operator $L_{\overline{\boldsymbol\gamma},x}^{w_1,w_2}$. For $u,\lambda\in i\R$ and $\rho,\sigma\in\R$, complex conjugation of \eqref{eqn:gamma-parameter-correspondence} shows that $\overline{\boldsymbol\gamma}$ corresponds to the parameters $(\lambda,u,-\rho,-\sigma)$. Since the integrand in \eqref{eqn:psi-integral} depends on $\rho$ and $\sigma$ only through $\pm\rho$ and $\pm\sigma$, we have $\psi_{\overline{\boldsymbol\gamma}}=\psi_{\boldsymbol\gamma}$. Hence \eqref{eqn:psi-eigenvalue-transform-section} also holds with $\boldsymbol\gamma$ replaced by $\overline{\boldsymbol\gamma}$, with eigenvalue $v(w_1,w_2,\overline{\boldsymbol\gamma};y')=\overline{v(w_1,w_2,\boldsymbol\gamma;y')}$ for real $y'$. Consequently, for $g=\psi_{\boldsymbol\gamma}^{y'}$,
\[
\langle f,L_{\overline{\boldsymbol\gamma},x}^{w_1,w_2}g\rangle_{\boldsymbol\gamma,M}
=\overline{v(w_1,w_2,\overline{\boldsymbol\gamma};y')}\,\langle f,g\rangle_{\boldsymbol\gamma,M}
=v(w_1,w_2,\boldsymbol\gamma;y')\,\langle f,g\rangle_{\boldsymbol\gamma,M}.
\]
In the conjugate case the same computation holds with $L_{\overline{\boldsymbol\gamma},x}^{w_2,w_1}$ and $v(w_2,w_1,\overline{\boldsymbol\gamma};y')=\overline{v(w_1,w_2,\boldsymbol\gamma;y')}$, using the equation with $w_1$ and $w_2$ interchanged.

Applying the lemma to
$f=\psi_{\boldsymbol\gamma}^{y}$ and
$g=\psi_{\boldsymbol\gamma}^{y'}$ gives, in the positive case,
\begin{equation}
	\begin{split}
		&
		\Bigl(
		v(w_1,w_2,\boldsymbol\gamma;y)
		-
		v(w_1,w_2,\boldsymbol\gamma;y')
		\Bigr)
		\left\langle
		\psi_{\boldsymbol\gamma}^{y},
		\psi_{\boldsymbol\gamma}^{y'}
		\right\rangle_{\boldsymbol\gamma,M}
		=
		2		\left[
		\psi_{\boldsymbol\gamma}^{y},
		\psi_{\boldsymbol\gamma}^{y'}
		\right]_{\boldsymbol\gamma}(M).
	\end{split}
	\label{eqn:Green-identity-psi}
\end{equation}
The analogous identity in the conjugate case uses
$\langle\cdot,\cdot\rangle_{\boldsymbol\gamma,M,\epsilon}$ and is followed by
$\epsilon\downarrow0$.

\subsection{Weak orthogonality and the weighted isometry}

For functions on $(0,\infty)$, set
\[
\langle f,g\rangle_{\boldsymbol\gamma,M}^{+}
:=
\int_0^M f(x)\overline{g(x)}W_{\boldsymbol\gamma}(x)\dd x.
\]

\begin{proposition}
	\label{prop:weak-orthogonality-psi}
	Suppose that
	\begin{equation}
		u,\lambda\in i(w_0-w,w),
		\qquad
		\operatorname{Im}(\lambda+u)>w_0-w,
		\qquad
		\rho,\sigma\in\R.
		\label{eqn:transform-parameter-domain}
	\end{equation}
	Then, for every $f\in C_{\mathrm c}^\infty(0,\infty)$ and $y'>0$,
	\begin{equation}
		\lim_{M\to\infty}
		\int_0^\infty
		f(y)
		\left\langle
		\psi_{\boldsymbol\gamma}^{y},
		\psi_{\boldsymbol\gamma}^{y'}
		\right\rangle_{\boldsymbol\gamma,M}^{+}
		\dd y
		=
		\frac{f(y')}{W_{\widehat{\boldsymbol\gamma}}(y')}.
	\end{equation}
	In the conjugate case, the truncated pairing is understood as the limit of
	the corresponding indented-contour expression as $\epsilon\downarrow0$.
\end{proposition}

\begin{proof}
	Substitute the asymptotic expansion from
	Theorem~\ref{thm:psi-asymptotics-article} into the Casorati determinant in
	\eqref{eqn:Green-identity-psi}. The resulting expression is a finite sum of
	terms involving
	\[
	\frac{\sin(\alpha M(y-y'))}{y-y'}
	\qquad\text{and}\qquad
	\frac{\sin(\alpha M(y+y'))}{y+y'},
	\]
	together with a remainder tending to zero locally uniformly. Here we use that, by \eqref{eqn:psi-asymptotic-estimate}, the error terms are smaller than the leading terms by a factor $e^{-\alpha\varepsilon\operatorname{Re}x}$ with $\varepsilon>0$, so that their contribution to the Casorati bracket, which is of the order of the product of the two kernels times $W_{\boldsymbol\gamma}\sim e^{2\alpha(w-\operatorname{Im}(\lambda+u))x}$, tends to zero. Note that the assumptions \eqref{eqn:transform-parameter-domain} imply the hypotheses of Theorem~\ref{thm:psi-asymptotics-article}. The terms
	involving $y+y'$ converge weakly to zero on $(0,\infty)$, whereas those
	involving $y-y'$ form an approximate identity at $y=y'$. The coefficient of
	the limiting Dirac mass simplifies, using \eqref{eqn:transform-weight}, to
	$1/W_{\widehat{\boldsymbol\gamma}}(y')$.
	The normalization here is for the half-line pairing. Indeed, the leading asymptotics give
	$W_{\boldsymbol\gamma}(x)\sim
	(w_1w_2)^{-1}e^{2\alpha(w+iu+i\lambda)x}$.
	The two conjugate plane waves therefore contribute 
	\[
	\frac{2\pi}{\alpha}\,
	c'(\widehat{\boldsymbol\gamma};y')
	c'(\widehat{\boldsymbol\gamma};-y')
	=
	\frac{1}{W_{\widehat{\boldsymbol\gamma}}(y')},
	\]
	since $\alpha=2\pi/w_1w_2$. The full-line pairing has twice this coefficient.
\end{proof}

\begin{theorem}
	\label{thm:weighted-R-transform-isometry}
	Under the assumptions \eqref{eqn:transform-parameter-domain},
	$\mathcal R_{\widehat{\boldsymbol\gamma}}$ extends uniquely to a unitary map
	\[
	\mathcal R_{\widehat{\boldsymbol\gamma}}:
	\mathcal H_{\widehat{\boldsymbol\gamma}}
	\longrightarrow
	\mathcal H_{\boldsymbol\gamma},
	\qquad
	\mathcal R_{\widehat{\boldsymbol\gamma}}^{-1}
	=\mathcal R_{\widehat{\boldsymbol\gamma}}^*
	=\mathcal R_{\boldsymbol\gamma}.
	\]
\end{theorem}
\begin{proof}
	For $f,g\in C_{\mathrm c}^\infty(0,\infty)$, insert
	\eqref{eqn:weighted-R-transform} into
	$\langle\mathcal R_{\widehat{\boldsymbol\gamma}}f,
	\mathcal R_{\widehat{\boldsymbol\gamma}}g\rangle_{\boldsymbol\gamma}$.
	Fubini's theorem, the duality
	$\psi_{\widehat{\boldsymbol\gamma}}(x,y)=
	\psi_{\boldsymbol\gamma}(y,x)$, and
	Proposition~\ref{prop:weak-orthogonality-psi} give
	$\langle\mathcal R_{\widehat{\boldsymbol\gamma}}f,
	\mathcal R_{\widehat{\boldsymbol\gamma}}g\rangle_{\boldsymbol\gamma}
	=\langle f,g\rangle_{\widehat{\boldsymbol\gamma}}$, so
	$\mathcal R_{\widehat{\boldsymbol\gamma}}$ extends to an isometry by density.

	Since the kernel is real and satisfies the duality
	relation, integration against compactly supported test functions shows that
	$\mathcal R_{\boldsymbol\gamma}=\mathcal R_{\widehat{\boldsymbol\gamma}}^*$
	on $C_{\mathrm c}^\infty(0,\infty)$. Now we make use of the
	$W(D_4)$-symmetry of the kernel $\mathcal E$. 
	The definitions \eqref{eqn:E-psi-relation} and \eqref{eqn:c-prime} give
	\[
	\mathcal E(\boldsymbol\gamma;x,y)
	=
	\frac{\chi(\boldsymbol\gamma)}{\sqrt{w_1w_2}}\,
	\frac{\psi_{\boldsymbol\gamma}(x,y)}
	{c'(\boldsymbol\gamma;x)c'(\widehat{\boldsymbol\gamma};y)}.
	\]
	Since $-\mathrm{id}\in W(D_4)$ and
	$\chi(-\boldsymbol\gamma)=\chi(\boldsymbol\gamma)$, 
	\eqref{eqn:E-D4-symmetry} implies
	\begin{equation}
	\psi_{\boldsymbol\gamma}(x,y)
	=
	h_{\boldsymbol\gamma}(x)
	h_{\widehat{\boldsymbol\gamma}}(y)
	\psi_{-\boldsymbol\gamma}(x,y),
	\qquad
	h_{\boldsymbol\gamma}(x)
	:=
	\frac{c'(\boldsymbol\gamma;x)}
	{c'(-\boldsymbol\gamma;x)}.
	\label{eqn:weighted-sign-reversal-kernel}
	\end{equation}
	By the reflection identity for $G$,
	\[
	h_{\boldsymbol\gamma}(x)
	=
	G(\lambda\pm\rho\pm x)\,G(u\pm\sigma\pm x)>0,
	\qquad x>0,
	\]
	where positivity follows as for $W_{\boldsymbol\gamma}$, since
	$|\operatorname{Im}\lambda|,|\operatorname{Im}u|<w$; the same holds for
	$h_{\widehat{\boldsymbol\gamma}}$, which has $(-u,-\lambda)$ in place of
	$(\lambda,u)$. Moreover $W_{-\boldsymbol\gamma}=h_{\boldsymbol\gamma}^2W_{\boldsymbol\gamma}$, so
	$U_{\boldsymbol\gamma}f:=h_{\boldsymbol\gamma}^{-1}f$ is a unitary map
	$\mathcal H_{\boldsymbol\gamma}\to\mathcal H_{-\boldsymbol\gamma}$, and
	similarly for $\widehat{\boldsymbol\gamma}$. By
	\eqref{eqn:weighted-sign-reversal-kernel},
	\[
	\mathcal R_{\boldsymbol\gamma}
	=
	U_{\widehat{\boldsymbol\gamma}}^{-1}
	\mathcal R_{-\boldsymbol\gamma}
	U_{\boldsymbol\gamma}
	\]
	on $C_{\mathrm c}^\infty(0,\infty)$. The parameter
	$-\widehat{\boldsymbol\gamma}$ corresponds under
	\eqref{eqn:gamma-parameter-correspondence} to
	$(\lambda,u,\rho,\sigma)\mapsto(u,\lambda,-\rho,-\sigma)$, which preserves
	\eqref{eqn:transform-parameter-domain}. The first part of the proof applied to
	$-\widehat{\boldsymbol\gamma}$ therefore shows that
	$\mathcal R_{-\boldsymbol\gamma}$ is isometric, and hence so is
	$\mathcal R_{\boldsymbol\gamma}$.

	Thus $\mathcal R_{\widehat{\boldsymbol\gamma}}$ is an isometry whose adjoint
	is also an isometry, i.e.\ it is unitary, and
	$\mathcal R_{\widehat{\boldsymbol\gamma}}^{-1}
	=\mathcal R_{\widehat{\boldsymbol\gamma}}^*
	=\mathcal R_{\boldsymbol\gamma}$.
\end{proof}

\begin{remark}
Since $G(cw_1,cw_2;cz)=G(w_1,w_2;z)$ for $c\in\C\setminus\{0\}$ with
$(w_1,w_2),(cw_1,cw_2)\in\mathcal P$, rescaling all variables and parameters
by $c$ shows that the results of this section also hold for periods
$(cw_1,cw_2)$ with $w_1,w_2$ positive or a complex conjugate pair, i.e.\
whenever $w_1/w_2\in(0,\infty)$ or $|w_1/w_2|=1$. In that case the variables
range over the ray $c\,(0,\infty)$ and the parameter conditions are rotated
accordingly.
\end{remark}

\subsection{A two-sided normalized overlap transform}

We introduce a two-sided normalization of the overlap kernel for which the
input weight depends only on $(\rho,\lambda)$ and the output weight depends
only on $(\sigma,-\lambda)$. This compatibility is essential for the
multivariate extension in the next subsection.

Let $\boldsymbol\gamma$ be related to
$(\lambda,u,\rho,\sigma)$ by
\eqref{eqn:gamma-parameter-correspondence}, and set
\begin{equation}
	a_{u,\sigma}(x):=G(u+\sigma\pm x),
	\qquad
	b_{u,\rho}(y):=G(-u+\rho\pm y).
\end{equation}
For $u,\lambda\in i\R$, $\rho,\sigma\in\R$, and $x,y\in\R$, the conjugation
formula \eqref{eqn:G-symmetries} gives
\begin{equation}
	\overline{a_{u,\sigma}(x)}=G(u-\sigma\pm x),
	\qquad
	\overline{b_{u,\rho}(y)}=G(-u-\rho\pm y).
\end{equation}

Recall the weight $W_{\boldsymbol\gamma}$ from
\eqref{eqn:transform-weight}. We define
\begin{equation}
	\begin{split}
		V_{\rho}^{\lambda}(x)
		&:={}|a_{u,\sigma}(x)|^2W_{\boldsymbol\gamma}(x)
		=\frac{1}{w_1w_2}
		\frac{G(iw\pm2x)}{G(\lambda\pm\rho\pm x)}
	\end{split}
	\label{eqn:two-sided-input-weight}
\end{equation}
and
\begin{equation}
	\begin{split}
		V_{\sigma}^{-\lambda}(y)
		&:={}|b_{u,\rho}(y)|^2
		W_{\widehat{\boldsymbol\gamma}}(y)
		=\frac{1}{w_1w_2}
		\frac{G(iw\pm2y)}{G(-\lambda\pm\sigma\pm y)}.
	\end{split}
	\label{eqn:two-sided-output-weight}
\end{equation}
The second identities follow by applying
\eqref{eqn:transform-weight} to $\boldsymbol\gamma$ and
$\widehat{\boldsymbol\gamma}$, respectively. In particular, the first weight
is independent of $u,\sigma$, whereas the second is independent of
$u,\rho$.

Define the normalized overlap kernel by
\begin{equation}
	K_{\rho,\sigma}^{\lambda,u}(x,y)
	:=
	\frac{\psi_{\rho,\sigma}^{\lambda,u}(x,y)}
	{G(u+\sigma\pm x)G(-u+\rho\pm y)}.
	\label{eqn:two-sided-overlap-kernel}
\end{equation}
For $f\in C_{\mathrm c}^{\infty}(0,\infty)$, set
\begin{equation}
	(\mathscr F_{\rho,\sigma}^{\lambda,u}f)(y)
	:=
	\int_0^\infty
	f(x)K_{\rho,\sigma}^{\lambda,u}(x,y)
	V_{\rho}^{\lambda}(x)\dd x.
	\label{eqn:two-sided-overlap-transform}
\end{equation}
Thus, if $M_h$ denotes multiplication by $h$, then
\begin{equation}
	\mathscr F_{\rho,\sigma}^{\lambda,u}
	=
	M_{b_{u,\rho}^{-1}}
	\mathcal R_{\boldsymbol\gamma}
	M_{\overline{a_{u,\sigma}}}.
	\label{eqn:two-sided-overlap-factorization}
\end{equation}

\begin{theorem}
	\label{thm:two-sided-overlap-transform}
	Suppose that \eqref{eqn:transform-parameter-domain} holds. Then
	\eqref{eqn:two-sided-overlap-transform} extends uniquely to a unitary map
	\begin{equation}
		\mathscr F_{\rho,\sigma}^{\lambda,u}:
		L^2\bigl((0,\infty),V_{\rho}^{\lambda}(x)\dd x\bigr)
		\longrightarrow
		L^2\bigl((0,\infty),V_{\sigma}^{-\lambda}(y)\dd y\bigr).
		\label{eqn:two-sided-overlap-unitary-map}
	\end{equation}
	Its inverse is
	\begin{equation}
		((\mathscr F_{\rho,\sigma}^{\lambda,u})^{-1}g)(x)
		=
		\int_0^\infty
		g(y)\overline{K_{\rho,\sigma}^{\lambda,u}(x,y)}
		V_{\sigma}^{-\lambda}(y)\dd y,
		\label{eqn:two-sided-overlap-inverse}
	\end{equation}
	where for a general $g \in  L^2((0,\infty), V_\sigma^{-\lambda}(y)\dd y)$ the integral is understood as the
	$L^2$ limit of integrals truncated to $(0,R)$.
	The following Plancherel identity holds
	\begin{equation}
		\begin{split}
			&\int_0^\infty
			|(\mathscr F_{\rho,\sigma}^{\lambda,u}f)(y)|^2
			V_{\sigma}^{-\lambda}(y)\dd y
			=
			\int_0^\infty |f(x)|^2V_{\rho}^{\lambda}(x)\dd x.
		\end{split}
		\label{eqn:two-sided-overlap-Plancherel}
	\end{equation}
\end{theorem}

\begin{proof}
	By \eqref{eqn:two-sided-input-weight}, multiplication by
	$\overline{a_{u,\sigma}}$ is unitary from the input space in
	\eqref{eqn:two-sided-overlap-unitary-map} onto
	$\mathcal H_{\boldsymbol\gamma}$. Similarly,
	\eqref{eqn:two-sided-output-weight} shows that multiplication by
	$b_{u,\rho}^{-1}$ is unitary from
	$\mathcal H_{\widehat{\boldsymbol\gamma}}$ onto the output space. The
	factorization \eqref{eqn:two-sided-overlap-factorization} and
	Theorem~\ref{thm:weighted-R-transform-isometry} prove the asserted unitarity
	and \eqref{eqn:two-sided-overlap-Plancherel}. Taking inverses in
	\eqref{eqn:two-sided-overlap-factorization} gives
	\[
	(\mathscr F_{\rho,\sigma}^{\lambda,u})^{-1}
	=
	M_{\overline{a_{u,\sigma}}^{-1}}
	\mathcal R_{\widehat{\boldsymbol\gamma}}
	M_{b_{u,\rho}},
	\]
	which yields \eqref{eqn:two-sided-overlap-inverse}. 
\end{proof}

\begin{remark}
	The relation \eqref{eqn:Rren-psi-relation} also expresses
	\eqref{eqn:two-sided-overlap-transform} as a transform with kernel
	$R_{\mathrm{ren}}$. The additional factors
	$G(\pm x+i\gamma_0)$ and
	$G(\pm y+i\widehat\gamma_0)$ are absorbed into multiplication operators. Hence the resulting $R_{\mathrm{ren}}$-transform
	is unitarily equivalent to the transform above and satisfies the corresponding
	inverse and Plancherel formulas. The overlap normalization
	\eqref{eqn:two-sided-overlap-kernel} is more convenient below because it
	iterates directly along the factors of the multivariate overlap coefficient
	\eqref{eqn:multivariate-psi}.
\end{remark}

\begin{remark}\label{rem:two-sided-kernel-symmetries}
 	The normalized kernel inherits several symmetries of the overlap
	coefficient. The symmetry
	\eqref{eqn:psi-switch} of the overlap coefficient gives
	\begin{equation*}
		K_{\rho,\sigma}^{\lambda,u}(x,y)
		=
		K_{\sigma,\rho}^{-\lambda,-u}(y,x).
	\end{equation*}
	Now assume the unitary parameter regime $x,y,\rho,\sigma\in\R$ and $\lambda,u\in i\R$. In this case the overlap coefficient is real-valued, and since it is also even in the boundary parameters $\rho,\sigma$, \eqref{eqn:G-symmetries} gives
	\begin{equation*}
		\overline{
			K_{\rho,\sigma}^{\lambda,u}(x,y)
		}
		=
		K_{-\rho,-\sigma}^{\lambda,u}(x,y).
	\end{equation*}
	Combining gives
	\begin{equation*}
		\overline{
			K_{\rho,\sigma}^{\lambda,u}(x,y)
		}
		=
		K_{-\sigma,-\rho}^{-\lambda,-u}(y,x).
	\end{equation*}
	The reduced weights $V_{\rho}^{\lambda}$ and $V_{\sigma}^{-\lambda}$ are also even in $\rho$ and $\sigma$.
	Consequently, the kernel of the inverse transform in \eqref{eqn:two-sided-overlap-inverse}, taken with respect to the measure $V_{\sigma}^{-\lambda}(y)\dd y$, is exactly the kernel of the transform \eqref{eqn:two-sided-overlap-transform} with parameters $(-\sigma,-\rho,-\lambda,-u)$, whose input and output weights are $V_{-\sigma}^{-\lambda}=V_{\sigma}^{-\lambda}$ and $V_{-\rho}^{\lambda}=V_{\rho}^{\lambda}$. Note, however, that the parameters $(-\lambda,-u)$ do not satisfy the assumptions of Theorem~\ref{thm:two-sided-overlap-transform} in the conjugate case (where $\lambda,u\in i(0,w)$), so this theorem cannot be applied to that transform directly. The following corollary shows that this is not necessary.
\end{remark}

\begin{corollary}
\label{cor:inverse-dual-transform}
Under the assumptions of Theorem~\ref{thm:two-sided-overlap-transform}, let $\mathscr F_{-\sigma,-\rho}^{-\lambda,-u}$ denote the integral operator defined by \eqref{eqn:two-sided-overlap-transform} on $C_{\mathrm c}^\infty(0,\infty)$ with parameters $(-\sigma,-\rho,-\lambda,-u)$ in place of $(\rho,\sigma,\lambda,u)$, i.e.\ with kernel $K_{-\sigma,-\rho}^{-\lambda,-u}$ and weight $V_{-\sigma}^{-\lambda}=V_\sigma^{-\lambda}$. Then $\mathscr F_{-\sigma,-\rho}^{-\lambda,-u}$ extends uniquely to a unitary map
\[
L^2\bigl((0,\infty),V_{\sigma}^{-\lambda}(y)\dd y\bigr)\longrightarrow L^2\bigl((0,\infty),V_{\rho}^{\lambda}(x)\dd x\bigr),
\]
and
\[
\left(\mathscr F_{\rho,\sigma}^{\lambda,u}\right)^{-1}=\mathscr F_{-\sigma,-\rho}^{-\lambda,-u}.
\]
\end{corollary}
In particular, the transform with the parameters $(-\sigma,-\rho,-\lambda,-u)$ is unitary also in the conjugate case, although these parameters lie outside the domain of Theorem~\ref{thm:two-sided-overlap-transform}.
\begin{proof}
By Remark~\ref{rem:two-sided-kernel-symmetries} (the kernel identity $\overline{K_{\rho,\sigma}^{\lambda,u}(x,y)}=K_{-\sigma,-\rho}^{-\lambda,-u}(y,x)$), for $g\in C_{\mathrm c}^\infty(0,\infty)$ the right-hand side of \eqref{eqn:two-sided-overlap-inverse} equals
\[
\int_0^\infty g(y)K_{-\sigma,-\rho}^{-\lambda,-u}(y,x)V_\sigma^{-\lambda}(y)\dd y=\bigl(\mathscr F_{-\sigma,-\rho}^{-\lambda,-u}g\bigr)(x).
\]
Hence $\mathscr F_{-\sigma,-\rho}^{-\lambda,-u}$ coincides on $C_{\mathrm c}^\infty(0,\infty)$ with the unitary operator $(\mathscr F_{\rho,\sigma}^{\lambda,u})^{-1}$, which proves both assertions.
\end{proof}

\subsection{The multivariate overlap transform}

We now apply the preceding construction to each factor of the multivariate
overlap coefficient. Let
$\boldsymbol x=(x_1,\ldots,x_N)$ and
$\boldsymbol y=(y_1,\ldots,y_N)$, with
$x_0:=\rho$ and $y_{N+1}:=\sigma$.
By \eqref{eqn:multivariate-psi}, the $j$th local overlap factor is
\[
\psi_{x_{j-1},y_{j+1}}^{\lambda_j,u_j}(x_j,y_j).
\]
Replacing $(\rho,\sigma,\lambda,u)$ in
\eqref{eqn:two-sided-overlap-kernel} by
$(x_{j-1},y_{j+1},\lambda_j,u_j)$ gives
\begin{equation}
	K_j(\boldsymbol x,\boldsymbol y)
	:=
	\frac{
		\psi_{x_{j-1},y_{j+1}}^{\lambda_j,u_j}(x_j,y_j)
	}{
		G(u_j+y_{j+1}\pm x_j)
		G(-u_j+x_{j-1}\pm y_j)
	}.
\end{equation}
Define the multivariate kernel
\begin{equation}
	\boldsymbol K_{\rho,\sigma}^{\boldsymbol\lambda,\boldsymbol u}
	(\boldsymbol x,\boldsymbol y)
	:=
	\prod_{j=1}^N K_j(\boldsymbol x,\boldsymbol y).
\end{equation}
Since the product of the numerators is precisely the multivariate overlap
coefficient from \eqref{eqn:multivariate-psi}, we have
\begin{equation}
	\begin{split}
		&\boldsymbol K_{\rho,\sigma}^{\boldsymbol\lambda,\boldsymbol u}
		(\boldsymbol x,\boldsymbol y)
		=
		\frac{
			\boldsymbol\psi_{\rho,\sigma}^{\boldsymbol\lambda,\boldsymbol u}
			(\boldsymbol x,\boldsymbol y)
		}{
			\prod_{j=1}^N
			G(u_j+y_{j+1}\pm x_j)
			G(-u_j+x_{j-1}\pm y_j)
		}.
	\end{split}
\end{equation}
Thus $\boldsymbol K$ is the multivariate overlap coefficient divided by the
two normalizing factors supplied by every one-variable transform.

The product of the local input weights in
\eqref{eqn:two-sided-input-weight} is
\begin{equation}
	\boldsymbol V_{\rho}^{\boldsymbol\lambda}(\boldsymbol x)
	:=
	(w_1w_2)^{-N}
	\prod_{j=1}^N
	\frac{G(iw\pm2x_j)}
	{G(\lambda_j\pm x_{j-1}\pm x_j)}.
	\label{eqn:multivariate-two-sided-input-weight}
\end{equation}
The product of the local output weights in
\eqref{eqn:two-sided-output-weight} is
\begin{equation}
	\boldsymbol V_{\sigma}^{-\boldsymbol\lambda}(\boldsymbol y)
	:=
	(w_1w_2)^{-N}
	\prod_{j=1}^N
	\frac{G(iw\pm2y_j)}
	{G(-\lambda_j\pm y_j\pm y_{j+1})}.
	\label{eqn:multivariate-two-sided-output-weight}
\end{equation}
The crucial point is that the first weight depends only on the $x$-variables,
whereas the second depends only on the $y$-variables. 

For a compactly supported function $f$ on $(0,\infty)^N$, define
\begin{equation}
	\begin{split}
		&(\boldsymbol{\mathscr F}_{\rho,\sigma}^{
			\boldsymbol\lambda,\boldsymbol u}f)(\boldsymbol y):=
		\int_{(0,\infty)^N}
		f(\boldsymbol x)
		\boldsymbol K_{\rho,\sigma}^{\boldsymbol\lambda,\boldsymbol u}
		(\boldsymbol x,\boldsymbol y)
		\boldsymbol V_{\rho}^{\boldsymbol\lambda}(\boldsymbol x)
		\dd\boldsymbol x.
	\end{split}
	\label{eqn:multivariate-two-sided-overlap-transform}
\end{equation}
Here the iterated integral is taken in the order
$\dd x_N\cdots\dd x_1$.

\begin{theorem}[Multivariate overlap transform]
	\label{thm:multivariate-two-sided-overlap-transform}
	Suppose that $\rho,\sigma\in\R$ and that, for each
	$j=1,\ldots,N$, the assumptions of
	Theorem~\ref{thm:two-sided-overlap-transform} hold for
	$(\lambda_j,u_j)$. Then
	\eqref{eqn:multivariate-two-sided-overlap-transform} extends uniquely to a
	unitary map
	\begin{equation}
		\begin{split}
			&\boldsymbol{\mathscr F}_{\rho,\sigma}^{
				\boldsymbol\lambda,\boldsymbol u}:
			L^2\left((0,\infty)^N,
			\boldsymbol V_{\rho}^{\boldsymbol\lambda}(\boldsymbol x)
			\dd\boldsymbol x\right)
			\longrightarrow
			L^2\left((0,\infty)^N,
			\boldsymbol V_{\sigma}^{-\boldsymbol\lambda}(\boldsymbol y)
			\dd\boldsymbol y\right).
		\end{split}
	\end{equation}
	Its inverse is
	\begin{equation}
		\begin{split}
			&((\boldsymbol{\mathscr F}_{\rho,\sigma}^{
				\boldsymbol\lambda,\boldsymbol u})^{-1}g)(\boldsymbol x)
			=
			\int_{(0,\infty)^N}
			g(\boldsymbol y)
			\overline{
				\boldsymbol K_{\rho,\sigma}^{\boldsymbol\lambda,\boldsymbol u}
				(\boldsymbol x,\boldsymbol y)}
			\boldsymbol V_{\sigma}^{-\boldsymbol\lambda}(\boldsymbol y)
			\dd\boldsymbol y,
		\end{split}
		\label{eqn:multivariate-two-sided-overlap-inverse}
	\end{equation}
	where the iterated integral is applied in the order $\dd y_1\cdots\dd y_N$.
	For general $g \in L^2$ these iterated integrals are understood as
	the successive $L^2$ extensions of the one-variable inverses. Its
	Plancherel identity is
	\begin{equation}
		\begin{split}
			&\int_{(0,\infty)^N}
			|(\boldsymbol{\mathscr F}_{\rho,\sigma}^{
				\boldsymbol\lambda,\boldsymbol u}f)(\boldsymbol y)|^2
			\boldsymbol V_{\sigma}^{-\boldsymbol\lambda}(\boldsymbol y)
			\dd\boldsymbol y
			=
			\int_{(0,\infty)^N}
			|f(\boldsymbol x)|^2
			\boldsymbol V_{\rho}^{\boldsymbol\lambda}(\boldsymbol x)
			\dd\boldsymbol x.
		\end{split}
		\label{eqn:multivariate-two-sided-overlap-Plancherel}
	\end{equation}
\end{theorem}

\begin{proof}
	Apply Theorem~\ref{thm:two-sided-overlap-transform} to the
	variable $x_N$, with left boundary $x_{N-1}$ and right boundary
	$y_{N+1}=\sigma$. Its input weight is
	$V_{x_{N-1}}^{\lambda_N}(x_N)$, while its output weight is
	$V_{\sigma}^{-\lambda_N}(y_N)$. In particular, the output weight no longer
	depends on $x_{N-1}$.
	The hypotheses of the one-variable theorem apply at every stage because
	the intermediate boundary variables are positive real numbers and the
	appropriate one-variable conditions are imposed for each $j$.
	
	The remaining integral has the same form with $N$ replaced by $N-1$.
	Applying the one-variable theorem successively in the order
	$x_N,x_{N-1},\ldots,x_1$ proves the norm identity
	\eqref{eqn:multivariate-two-sided-overlap-Plancherel}. At the $j$th step the
	input weight is $V_{x_{j-1}}^{\lambda_j}(x_j)$ and the output weight is
	$V_{y_{j+1}}^{-\lambda_j}(y_j)$, so the factors assemble exactly into
	\eqref{eqn:multivariate-two-sided-input-weight} and
	\eqref{eqn:multivariate-two-sided-output-weight}.

	Let us make the iteration precise. Let $f$ be continuous with compact support in $(0,\infty)^N$. After the transforms in $x_N,\ldots,x_{j+1}$ we have a function $f_j(x_1,\ldots,x_j;y_{j+1},\ldots,y_N)$, continuous in all its variables because the kernels are continuous, and still compactly supported in $(x_1,\ldots,x_j)$, since these variables have not been touched. The $j$th step applies the one-variable transform in $x_j$, with the real boundary parameters $x_{j-1}$ and $y_{j+1}$, to a function that is continuous with compact support in $x_j$. For such functions the integral formula holds pointwise, and Theorem~\ref{thm:two-sided-overlap-transform} applies, because $C_{\mathrm c}^\infty(0,\infty)$ is dense in $C_{\mathrm c}(0,\infty)$ with respect to uniform convergence with supports in a fixed compact set, which implies convergence in the $L^2$-norms of the weights and of the transformed functions. The remaining variables $(x_1,\ldots,x_{j-1},y_{j+1},\ldots,y_N)$ enter only as parameters, all functions are measurable in them, and the one-variable Plancherel identity can be integrated over these variables by Tonelli's theorem. In this way \eqref{eqn:multivariate-two-sided-overlap-Plancherel} holds on the dense class $C_{\mathrm c}((0,\infty)^N)$, and the extension to $L^2$ is by continuity.

	Applying the inverse one-variable transforms in the reverse order
	$y_1,y_2,\ldots,y_N$ gives
	\eqref{eqn:multivariate-two-sided-overlap-inverse} on the same dense
	subspace. The norm identity and inverse relation extend by continuity, proving
	that the resulting operator is unitary.
\end{proof}

\begin{remark} 
	Unlike the coproduct construction of the multivariate difference operators in Section~\ref{subsec:multivariate-difference-equations}, the transform construction does not require $u_1=\cdots=u_N$. This is because it is obtained by iterating one-variable transforms rather than by a coproduct-eigenfunction argument.
\end{remark}

\appendix
\section{The hyperbolic gamma function}
\label{sec:hyperbolic-gamma}

In this appendix we collect the properties of the hyperbolic gamma function, following \cite{RuijsenaarsADE,R1},
that are used in the paper. Throughout, $(w_1,w_2)\in\mathcal P$, and $\mathcal P$, $w$ and $\alpha$ are as in Section~\ref{subsec:notation}.

For $z$ in the strip $\{\left|\operatorname{Im}z\right|<\operatorname{Re}w \}$ define
\begin{equation*}
	g(w_1,w_2;z)
	:=
	\int_0^\infty
	\left(
	\frac{\sin(2yz)}
	{2\sinh(w_1y)\sinh(w_2y)}
	-
	\frac{z}{w_1w_2y}
	\right)
	\frac{\dd y}{y}.
\end{equation*}
The integral is locally uniformly convergent in this strip and defines an
analytic function of $z$. The hyperbolic gamma function is defined there by
\begin{equation*}
	G(w_1,w_2;z):=
	\exp\bigl(ig(w_1,w_2;z)\bigr).
\end{equation*}
It follows in particular that $G$ is analytic and nonzero on the given strip.
When the dependence on $w_1,w_2$ is clear, we write $G(z)$ instead of
$G(w_1,w_2;z)$.

\subsection{Symmetries and difference equations}

The definition immediately gives
\begin{equation}
	G(w_1,w_2;-z)
	=
	G(w_1,w_2;z)^{-1}
	=
	\overline{
		G(\overline{w_1},\overline{w_2};\overline z)
	}.
	\label{eqn:G-symmetries}
\end{equation}
Moreover,
\begin{equation*}
	G(w_1,w_2;z)=G(w_2,w_1;z).
\end{equation*}
In particular, if $w_1,w_2>0$ or $w_2=\overline{w_1}$, then
$\{\overline{w_1},\overline{w_2}\}=\{w_1,w_2\}$ and $w\in\R$, so that
\begin{equation}
	\overline{G(z)}=G(-\overline z)=G(\overline z)^{-1}.
	\label{eqn:G-conjugation}
\end{equation}

The hyperbolic gamma function satisfies the difference equations
\begin{equation}
	\frac{G(z+iw_1/2)}{G(z-iw_1/2)}
	=
	2\cosh\frac{\pi z}{w_2}
	\label{eqn:G-difference-w1}
\end{equation}
and
\begin{equation*}
	\frac{G(z+iw_2/2)}{G(z-iw_2/2)}
	=
	2\cosh\frac{\pi z}{w_1}.
\end{equation*}
These identities initially hold whenever both arguments belong to the
defining strip.

\subsection{Meromorphic continuation, zeros, and poles} \label{app:zeros-and-poles}

The difference equations extend $G$ to a meromorphic function on $\C$. In fact, $G$ is meromorphic jointly in
$(w_1,w_2,z)\in\mathcal P\times\C$ \cite[App.~A]{R1}, and the identities above hold as identities of meromorphic functions.

The hyperbolic gamma function is scale invariant:
\begin{equation}
	G(cw_1,cw_2;cz)=G(w_1,w_2;z)
	\label{eqn:G-scaling}
\end{equation}
for $c>0$ \cite[(A.6)]{R1}. By analytic continuation in $c$, this extends to all
$c\in\C\setminus\{0\}$ with $(cw_1,cw_2) \in \mathcal P$.

The zeros of $G$ are contained in
\begin{equation*}
	Z_+
	:=
	i\left\{
	w+kw_1+lw_2
	\mathrel{}\middle|\mathrel{}
	k,l\in\N\cup\{0\}
	\right\},
\end{equation*}
and its poles are contained in
\begin{equation*}
	Z_-:=-Z_+.
\end{equation*}
If $w_1/w_2\not\in\Q$, these zeros and poles are simple. In any case,
$G$ has a simple zero at $iw$ and a simple pole at $-iw$, since $w+kw_1+lw_2=w$ with $k,l\in\N\cup\{0\}$ forces $k=l=0$. Moreover, since $\operatorname{Re}(w+kw_1+lw_2)>\operatorname{Re}w$ unless $k=l=0$, the only zero of $G$ on the line $\operatorname{Im}z=\operatorname{Re}w$ is $iw$, and the only pole of $G$ on the line $\operatorname{Im}z=-\operatorname{Re}w$ is $-iw$.

The residue at the pole nearest the defining strip is
\begin{equation*}
	\Res_{z=-iw}G(z)
	=
	\frac{i\sqrt{w_1w_2}}{2\pi},
\end{equation*}
where $\sqrt{w_1w_2}$ denotes the principal branch; recall that
$\arg(w_1w_2)\in(-\pi,\pi)$ on $\mathcal P$. It is positive for positive $w_1,w_2$ and equals $|w_1|$ when $w_2=\overline{w_1}$. By
\eqref{eqn:G-symmetries}, this is equivalent to
\begin{equation*}
	G(z)
	=
	\frac{i\sqrt{w_1w_2}}{2\pi(z+iw)}
	+O(1)
	\qquad
	\text{as }z\longrightarrow-iw.
\end{equation*}
Similarly,
\begin{equation*}
	G(z)
	=
	\frac{2\pi i}{\sqrt{w_1w_2}}(z-iw)
	+O\bigl((z-iw)^2\bigr)
	\qquad
	\text{as }z\longrightarrow iw.
\end{equation*}

\subsection{Asymptotic estimates}
In this subsection we assume that $w_1,w_2>0$ or $w_2=\overline{w_1}$; then $\alpha>0$. Define
\[
w_0:=
\begin{cases}
	\min\{w_1,w_2\},& w_1,w_2>0,\\[0.2cm]
	w,& \overline{w_1}=w_2.
\end{cases}
\]
In the second case $w$ is real and positive. Fix
$\varsigma\in(0,1)$, and set
\begin{equation*}
	\widetilde g(z)
	:=
	-\frac{\alpha}{4}z^2
	-\frac{\alpha}{48}(w_1^2+w_2^2).
\end{equation*}
Let $K\subset\R$ be compact. Uniformly for
$\operatorname{Im}z\in K$, and away from fixed neighborhoods of the pole
and zero sequences, we have
\begin{equation*}
	G(z)
	=
	e^{i\widetilde g(z)}
	\left(
	1+O\left(e^{-\alpha\varsigma w_0\operatorname{Re}z}\right)
	\right)
	\qquad
	\text{as }\operatorname{Re}z\longrightarrow+\infty.
\end{equation*}
This is \cite[Thm.~A.1]{R1}, where the error is bounded by
$C\exp(-2\pi\sigma\operatorname{Re}z/a_m)$ with $\sigma\in(\frac12,1)$ and
$a_m=\max_j|w_j|^2/\operatorname{Re}w_j$; in both cases considered here
$2\pi/a_m=\alpha w_0$. (For $\varsigma\le\frac12$ the bound above is weaker.)
Using $G(-z)=G(z)^{-1}$, this gives
\begin{equation*}
	G(z)
	=
	e^{-i\widetilde g(-z)}
	\left(
	1+O\left(e^{-\alpha\varsigma w_0|\operatorname{Re}z|}\right)
	\right)
	\qquad
	\text{as }\operatorname{Re}z\longrightarrow-\infty.
\end{equation*}

In particular, if $a,b\in\C$ range over compact sets, then, uniformly for
$\operatorname{Im}z$ in compact sets,
\begin{equation*}
	\frac{G(z-a)}{G(z-b)}
	=
	\exp\left(
	-\frac{i\alpha}{2}z(b-a)
	+
	\frac{i\alpha}{4}(b^2-a^2)
	\right)
	\left(
	1+O\left(e^{-\alpha\varsigma w_0\operatorname{Re}z}\right)
	\right)
\end{equation*}
as $\operatorname{Re}z\to+\infty$.

The corresponding estimate for absolute values is
\begin{equation*}
	0<C_-
	<
	\left|
	\frac{G(z-a)}{G(z-b)}
	\right|
	\exp\left(
		\frac{\alpha}{2}
	\operatorname{Im}(a-b)
	|\operatorname{Re}z|
	\right)
	<C_+
\end{equation*}
for $|\operatorname{Re}z|$ sufficiently large and
$\operatorname{Im}z$ in a fixed compact set. The constants $C_-,C_+>0$
may be chosen uniformly when $a$ and $b$ range over compact sets and the
arguments remain a fixed positive distance from the zeros and poles of $G$.

A special case that is used repeatedly is
\begin{equation*}
	G(iy\pm z)
	=
	\begin{cases}
		e^{\alpha yz}
		\left(1+O\left(e^{-\alpha\varsigma w_0\operatorname{Re}z}\right)\right),
		&\text{as }\operatorname{Re}z\longrightarrow+\infty,\\
		e^{-\alpha yz}
		\left(1+O\left(e^{-\alpha\varsigma w_0|\operatorname{Re}z|}\right)\right),
		&\text{as }\operatorname{Re}z\longrightarrow-\infty,
	\end{cases}
\end{equation*}
uniformly for $y$ and $\operatorname{Im}z$ in compact sets. In
particular,
\begin{equation*}
	0<C_-
	<
	|G(iy\pm z)|
	e^{-\alpha\operatorname{Re}(y)|\operatorname{Re}z|}
	<C_+
\end{equation*}
for $|\operatorname{Re}z|$ sufficiently large.

More generally, if $a_1,\ldots,a_m,b_1,\ldots,b_n$ range over compact
sets, the estimates above may be applied factorwise to finite products of
the form
\[
\frac{\prod_{j=1}^mG(z-a_j)}
{\prod_{k=1}^nG(z-b_k)}.
\]
All resulting error estimates are locally uniform in the parameters, as
long as the arguments remain a fixed positive distance from the pole and
zero sequences.

\subsection{Relation with Faddeev's quantum dilogarithm}
\label{app:quantum-dilogarithm}
The hyperbolic gamma function is closely related to Faddeev's noncompact quantum dilogarithm; see \cite{FKV}. Since several normalizations of the quantum dilogarithm occur in the literature, we state the convention used here explicitly. Set $b:=w_1/\sqrt{w_1w_2}$, where we use the principal branch of the square root, thus $b^2=w_1/w_2$. Faddeev's quantum dilogarithm is defined by
\[ 
\Phi_b(x) := \exp\left( \frac14 \int_{\R+i0} \frac{e^{-2ixt}} {\sinh(bt)\sinh(t/b)} \frac{\dd t}{t} \right),
\qquad |\operatorname{Im}x|<\tfrac12\operatorname{Re}(b+b^{-1}),
\] 
and by meromorphic continuation elsewhere. Note that the scale invariance \eqref{eqn:G-scaling} with $c=1/\sqrt{w_1w_2}$, for which $(cw_1,cw_2)=(b,b^{-1})\in\mathcal P$, gives
\[
G(w_1,w_2;z)=G\bigl(b,b^{-1};z/\sqrt{w_1w_2}\bigr),
\]
so it suffices to compare $G(b,b^{-1};\cdot)$ with $\Phi_b$.
Then 
\begin{equation} 
	G(w_1,w_2;z) = \Phi_b(0) \frac{ \exp\left(\dfrac{i\pi z^2}{2w_1w_2}\right) }{ \Phi_b\left(\dfrac{z}{\sqrt{w_1w_2}}\right) }, \label{eqn:G-quantum-dilogarithm} 
\end{equation} 
where 
\begin{equation*} 
	\Phi_b(0) = \exp\left( \frac{i\pi}{24}\left(b^2+b^{-2}\right) \right). 
\end{equation*} 
The Gaussian factor in \eqref{eqn:G-quantum-dilogarithm} accounts for the difference between the difference equations of $G$ and those of $\Phi_b$.

\end{document}